\documentclass[pdflatex,sn-mathphys]{sn-jnl}
\usepackage{xcolor}

\usepackage{subfig}
\usepackage{algorithm}
\usepackage{algpseudocode}
\renewenvironment{table}[1][]%
{\tableorg[#1]%
\tablebodyfont%
\renewcommand\footnotetext[2][]{{\removelastskip\vskip3pt%
\let\tablebodyfont\tablefootnotefont%
\hskip0pt\if!##1!\else{\smash{$^{##1}$}}\fi##2\par}}%
}{\endtableorg}

\jyear{2022}%

\theoremstyle{thmstyleone}%
\theoremstyle{thmstyletwo}%

\theoremstyle{thmstylethree}%

\begin{document}

\title[A High-Order Hybrid-Spectral Incompressible Navier-Stokes Model For Nonlinear Water Waves]{A High-Order Hybrid-Spectral Incompressible Navier-Stokes Model For Nonlinear Water Waves}


\author*[1,2]{\fnm{Anders} \sur{Melander}}\email{adame@dtu.dk}
\author[1]{\fnm{Max E.} \sur{Bitsch}}
\author[1]{\fnm{Dong} \sur{Chen}}
\author[1,2]{\fnm{Allan P.} \sur{Engsig-Karup}}

\affil[1]{\orgdiv{Department of Applied Mathematics and Computer Science}, \orgname{Technical University of Denmark}, \orgaddress{\street{Richard Petersens Plads}, \city{Kgs. Lyngby}, \postcode{2800}, \state{State}, \country{Denmark}}}
\affil[2]{\orgdiv{Center for Energy Resources Engineering (CERE)}}


\abstract{
We present a new high-order accurate computational fluid dynamics model based on the incompressible Navier-Stokes equations with a free surface for the accurate simulation of nonlinear and dispersive water waves in the time domain. The spatial discretization is based on Chebyshev polynomials in the vertical direction and a Fourier basis in the horizontal direction, allowing for the use of the fast Chebyshev and Fourier transforms  for the efficient computation of spatial derivatives. The temporal discretization is done through a generalized low-storage explicit 4th order Runge-Kutta, and for the scheme to conserve mass and achieve high-order accuracy, a velocity-pressure coupling needs to be satisfied at all Runge-Kutta stages. This result in the emergence of a Poisson pressure problem that constitute a geometric conservation law for mass conservation. The occurring Poisson problem is proposed to be solved efficiently via an accelerated iterative solver based on a geometric $p$-multigrid scheme, which takes advantage of the high-order polynomial basis in the spatial discretization and hence distinguishes itself from conventional low-order numerical schemes. We present numerical experiments for validation of the scheme in the context of numerical wave tanks demonstrating that the $p$-multigrid accelerated numerical scheme can effectively solve the Poisson problem that constitute the computational bottleneck, that the model can achieve the desired spectral convergence, and is capable of simulating wave-propagation over non-flat bottoms with excellent agreement in comparison to experimental results. 
}




\maketitle

\section{Introduction}\label{sec1}
Modelling of water waves are important across applications in offshore engineering, e.g. for estimation of sea states and wave-induced loads on structures derived from accurate kinematics and assessment of impact of environmental conditions etc. 
It is therefore of interest to be able to accurately and robustly perform numerical simulations of nonlinear and dispersive wave transformation and propagation subject to the influence of seabeds. A classical way to do this has been to use models such as Boussinesq-type equations \cite{Madsen91,Wei95,Li99}. These modelling methods reduce the dimension of the computational domain, therefore allowing for much faster computation compared to solving the more complete Navier-Stokes equations that take into account viscosity and irrotational effects. However, any such modeling provides tradeoffs between modeling fidelity in terms of physics and speed of computations \cite{Brocchini2013}. For example, Boussinesq-type modelling can handle nonlinear and dispersive waves (e.g. see \cite{ESKILSSON2014261}), however, it is generally only applicable to shallow and intermediate water regions in coastal areas and comes with a limitation that it is not straightforward to include wave-structure interaction despite recent years progress, e.g. see \cite{BOSI2019,beck2023numerical}. This is therefore an issue for handling floating bodies in connection with offshore engineering applications. As computational power has increased, so has the interest in using the more complete non-hydrostatic models which better take into account the depth variation relevant for sea state estimations and accurately estimating kinematics in the entire water column to produce improved load estimates in areas of interest for offshore structures. These include models based on potential flow theory, such as the fully nonlinear potential flow (FNPF) model. In the FNPF model a Laplace problem is solved to resolve the scalar velocity potential in the full domain, and it has been shown to work well for long time wave simulations and at regional scales \cite{Dommermuth_Yue_1987,Ducrozet07, Glimberg13, Allan09}. Moreover, FNPF modelling has also been shown to handle the complex geometry of offshore structures and wave-structure interaction \cite{Ducrozet12}. However, the FNPF model brings with it its own set of limitations, as the model is only valid for irrotational flow, and does not take viscosity into account \cite{Allan2}. 
On the other hand, direct Navier-Stokes equations (NSE) are widely and routinely used in applied computational fluid dynamics (CFD) communities and constitute a more complete model over the simplified models and allow for modeling rotational and viscous flow effects, turbulence, boundary layer separation (fx. for bluff bodies), wave-structure interaction, etc.  However, this comes at the cost of higher computational time in comparison to the earlier mentioned methods. It is noted that to lower the cost of CFD it is possible to utilize more efficient wave propagation models such as FNPF models as far-field solvers \cite{Allan09,Dommermuth_Yue_1987} and coupling these to conventional CFD solvers, e.g. see \cite{PAULSEN2014,Verbruggh2016,CHOI2023100510}. Solving the full set of NSE directly with a free surface is conventionally also done with numerical schemes such as Volume of Fluid (VOF) \cite{Hirt81} and level-set method \cite{Osher88, GroossHesthaven06} that implicitly capture the free surface waves. However, VOF is inherently a low order numerical scheme, which therefore can lead to a smearing of the free-surface due to excessive numerical diffusion and hence it is not very suitable for long-time integrations as can be concluded from analysis of numerical phase errors for wave propagation schemes \cite{KO72}. Furthermore, such conventional numerical schemes for many applications induce high computational costs due to spatio-temporal resolution requirements. An alternative method is to represent the free surface as a kinematic boundary condition to be able to track it explicitly, which also allows for the use of a $\sigma$-transform to map to a time-constant domain, which simplifies computations and reduces computational costs. This $\sigma$-transformation method seem to have been first proposed by Phillips \cite{Phillips57}. Later both Li \& Flemming \cite{Li01} and Chern et al. \cite{Chern99} applied this transformation to water wave modeling. Li \& Flemming used the method along with a projection method similar to Chorin \cite{Chorin68} to resolve the pressure through a Poisson boundary value problem. Engsig-Karup et al \cite{Allan24} designed a high-order Finite-Difference scheme to solve the free surface Navier-Stokes problem using high-order time-stepping with a mixed-stage Poisson problem defined to ensure a divergence free flow. Pan et al \cite{Pan21} proposed a discretization based on the discontinuous Galerkin finite element method (DGFEM), which allowed for a high-order accurate numerical scheme. The use of high-order, spectral accurate methods have shown good results for models, such as FNPF \cite{Christiansen13,Klahn20,Allan16,Bonnefoy10} and Boussinesq-type modelling \cite{MADSEN_BINGHAM_LIU_2002,Allan06}. These results show that having access to high-order spectral accurate solvers is beneficial for efficient solutions at large scale, as the high-order polynomial basis has a high-cost effectiveness. The high-order accuracy also lead to minimal numerical diffusion. This is especially important for time-dependent problems due to the large cost associated with solving equations at every time step and the accumulated effects of numerical diffusion errors. Spectral methods have also been shown to be effective in solving the incompressible Navier-Stokes \cite{Ku87,Spalart91}, however, little research has been done with these methods when including a free surface for water wave simulations. Applying these spectral methods to free surface CFD modelling of wave propagation is therefore scientifically relevant to benefit from the efficiency and cost effectiveness of spectral methods that have already been demonstrated for the related FNPF-based models, e.g. see  \cite{Dommermuth_Yue_1987,AGNON1999527,CLAMOND_GRUE_2001,Christiansen13}.

For CFD solvers, the Poisson problem for pressure generally presents a computational bottleneck and has the largest impact on runtime performance, and therefore efficient iterative solvers are required for solving efficiently the large linear algebraic systems that stem from the numerical discretization. For example, iterative multigrid solvers are among the most efficient and scalable algorithms and exploit the fact that classic stationary iterative methods are great at reducing high-frequency errors \cite{Trottenberg01}. Therefore, by designing solution procedure that utilize solving systems of algebraic equations on a hierarchy of multiple grids, this allows one to maintain the effectiveness of classic stationary iterative methods due to low-frequency errors appearing as high-frequency errors on coarser grids. Geometric multigrid methods have been shown to be effective for water wave problems, e.g. for free surface NSE modelling \cite{Li01} or FNPF modelling \cite{Allan09}. Moreover, the use of high-order accurate numerical methods allows for the use of geometric $p$-multigrid methods, which extends the multigrid method to exploit a hierarchy of different polynomial orders. Engsig-Karup \& Laskowski \cite{AllanLaskowski} showed good results using $p$-multigrid as preconditioning for the stationary defect correction (PDC) method \cite{EngsigKarup2014} and preconditioned conjugate gradient
(PCG) methods in FNPF models \cite{EngsigKarup2014,Allan16}. Such $p$-multigrid methods have seen limited exploration in the context of free surface Navier-stokes equations and is therefore explored in this work. 

\subsection{Contributions}

In this work, our aim is two-fold. Design a new high-order spectrally accurate numerical model for solving the incompressible Navier-Stokes equations with a free surface by combining a high-order hybrid-spectral spatial discretization with a $p$-multigrid accelerated generalized minimal residual method (GMRES) iterative scheme to efficiently solve the occurring Poisson problem for pressure. The equations are solved on a $\sigma$-transformed time-constant computational grid to reduce computational costs. To advance the equations in time a general $s$-stage low-storage explicit Runge-Kutta method (LSERK) is used and designed to satisfy a pressure-velocity coupling in line with the approach proposed by Chorin \cite{Chorin68} that ensures mass conservation and is suitable for high-order accurate temporal integration. Through convergence studies we demonstrate spectral convergence histories for shallow, intermediate and deep water wave simulation of nonlinear stream function waves. These results may be useful for benchmarking related numerical solvers.

\section{Mathematical formulation}

\begin{figure}[H]
    \centering
    \hspace*{-1cm}
    \includegraphics[width=100mm]{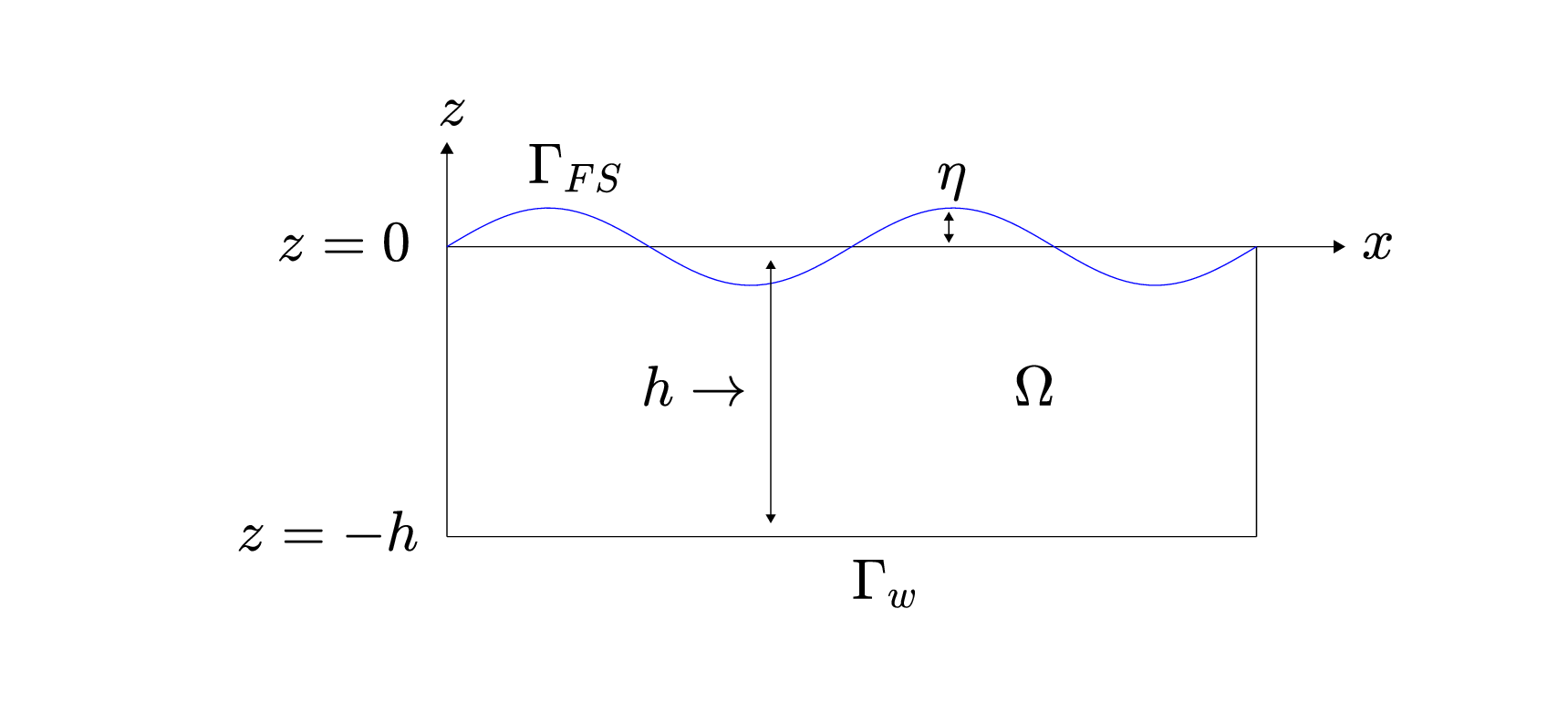}
    \caption{Illustration of a two-dimensional fluid domain, $\Omega$, and associated boundaries, $\Gamma_w$ and moving boundary $\Gamma_{FS}$. The free surface is denoted by $\eta(x, t)$ and the water depth by $h(x)$.}
    \label{fig:domain}
\end{figure}

The motion of a viscous fluid is governed by the Navier-Stokes equations in a spatial domain $\Omega$ and temporal domain $T$. To model the evolution of water wave motion a free surface kinematic boundary condition is included in the governing equations. Thereby, the governing equations in the 3D space of the fluid volume consists of the dynamic conditions describing changes in the velocity field, a kinematic free surface boundary condition, an additional mass continuity equation for incompressible flow, and momentum equations 
\begin{subequations}
    \begin{align}
        & \nabla \cdot \textbf{u} = 0, & \texttt{in } \Omega, \\ 
        &\frac{\partial \textbf{u}}{\partial t} + \textbf{u} \cdot\nabla  \textbf{u} = - \frac{1}{\rho} \nabla p + \textbf{g} + \nu \nabla^2 \textbf{u}, & \texttt{in } \Omega \times T,\label{eq:moms}
    \end{align} \label{eq:NSE}%
\end{subequations}
where $\textbf{u}$ are the velocities $[u,v,w]^T$ $\left[\frac{m}{s}\right]$, $\nabla$ is the gradient operator $[\frac{\partial}{\partial x},\frac{\partial}{\partial y},\frac{\partial}{\partial z}]^T$, $\rho$ is the density $\left[\frac{kg}{m^3}\right]$, $p$ is the pressure $\left[\text{Pa}\right]$, $\textbf{g}$ are the external forces $[0,0,-g]^T$, with $g$ being the gravitational acceleration $\left[\frac{m}{s^2}\right]$ and $\nu$ is the kinematic viscosity $\left[\frac{m^2}{s}\right]$. In this paper, we will use $\textbf{g} = [0,0,9.81]^T \left[\frac{m^2}{s}\right]$ for the gravitational acceleration and $\nu = 10^{-6}\left[\frac{m^2}{s}\right]$ for the kinematic viscosity.\\
We are interested in computing the NSE within some fluid domain $\Omega$ for simulating water waves. The evolution of the free surface $\eta(x,y,t)$ is governed by a free surface kinematic boundary condition
\begin{equation}
    \frac{\partial\eta}{\partial t} = -\tilde{u}\frac{\partial \eta}{\partial x}-\tilde{v}\frac{\partial \eta}{\partial y}+\tilde{w}, \quad \texttt{on } \Gamma_{FS}.
    \label{eq:kinematicFS}
\end{equation}
Here "$\sim$" denotes the solution at the free surface $z=\eta$, e.g. $\tilde{u} = u(x,y,\eta)$, see Figure \ref{fig:domain}.\\
Moreover, we make the assumption that the pressure is zero at the free surface, $p=0$ Pa, and choose to split the pressure into hydrostatic and hydrodynamic parts, such that 
\begin{equation}
    p = p_D + p_S, \quad p_S = \rho g(\eta-z), 
\end{equation}
where the static pressure, $p_S$, must be zero at the surface due to the linear dependence on the depth, and hence $p_D$ must also be zero at the surface to fulfill the earlier assumption about the total pressure at the free surface. 
The domain is bounded by impermeability conditions at walls, stated as 
\begin{equation}
    \textbf{n} \cdot \textbf{u} = 0 \quad \text{on} \quad \Gamma_w,
\end{equation}
where $\textbf{n}$ is the outward pointing normal at domain boundary points in question and $\Gamma_w$ represents the solid boundaries.

Furthermore, we note that if we take the divergence $\nabla\cdot$ of the momentum equations \eqref{eq:moms} and assume commuting operators, we obtain
\begin{align}
        \frac{\partial (\nabla\cdot \textbf{u})}{\partial t} + \nabla\cdot (\textbf{u} \cdot\nabla  \textbf{u}) = - \frac{1}{\rho} \nabla^2 p + \nabla\cdot \textbf{g} + \nu \nabla\cdot(\nabla^2 \textbf{u}),\label{eq:moms2}
\end{align} 
which can be interpreted as a geometric conservation law. The continuous formulation given in this section describe a divergence free flow for all times, and therefore we assume
\begin{align}
\frac{\partial (\nabla\cdot \textbf{u})}{\partial t} = 0,
\label{eq:velocitypressurecont}
\end{align}
resulting in a velocity-pressure coupling defined by (\ref{eq:moms2}).

\subsection{A $\sigma$-Coordinate Transformation} \label{sec:CT}
The spatial domain $\Omega(t)$ is bounded from above by the free surface $\eta(x,y,t)$, hence the fluid domain may vary over time. It is therefore useful to introduce a coordinate transform for the vertical coordinate in the form of a $\sigma$-transformation given as 
\begin{equation}
    \sigma = \frac{z+h(x,y)}{d(x,y,t)}, \quad 0 \leq \sigma \leq 1, \label{eq:sigtrans}
\end{equation}
where $d(x,y,t) = h(x,y)+\eta(x,y,t)$ represents the height of the water column and $h(x,t)$ is the depth of the seabed measured from the still-water level at $z=0$ m. We let
\begin{equation}
    \tau = t, \quad \chi = x, \quad \gamma = y , \quad \sigma = \sigma(t,x,y,z) ,
\end{equation}
be the transformed coordinates, where only the vertical coordinate is changed from the original domain.

To transform the NSE to the $\sigma$-domain, we first let
\begin{equation}
    f(t,x,y,z) = f(\tau,\chi,\gamma,\sigma(t,x,y,z) ),
\end{equation}
represent any of the equations in the full NSE formulation. Applying the chain rule gives us the transformed physical derivatives as
\begin{equation}
    \begin{aligned}
        & \frac{\partial f}{\partial t} = \frac{\partial f}{\partial \tau} + \frac{\partial f}{\partial \sigma} \frac{\partial \sigma}{\partial t}, \quad 
        \frac{\partial f}{\partial x} = \frac{\partial f}{\partial \chi} + \frac{\partial f}{\partial \sigma} \frac{\partial \sigma}{\partial x},\\ 
        & \frac{\partial f}{\partial y} = \frac{\partial f}{\partial \gamma} + \frac{\partial f}{\partial \sigma} \frac{\partial \sigma}{\partial y}, \quad
        \frac{\partial f}{\partial z} = \frac{\partial f}{\partial \sigma} \frac{\partial \sigma}{\partial z}.
    \end{aligned} \label{eq:1_chain}
\end{equation}
The $\sigma$-transformed derivatives can now be substituted into the Navier-Stokes formulations to find the $\sigma$-transformed equations. Note that we also include the pressure splitting here. We get
\begin{subequations} 
    \begin{align}
            & \frac{\partial u}{\partial \chi} + \frac{\partial u}{\partial \sigma}\frac{\partial \sigma}{\partial x} + \frac{\partial v}{\partial \gamma} +\frac{\partial v}{\partial \sigma}\frac{\partial \sigma}{\partial y} + \frac{\partial w}{\partial \sigma}\frac{\partial \sigma}{\partial z} = 0, \\ \label{eq:1_STNV_1}
            &  \frac{\partial u}{\partial \tau} + u\frac{\partial u}{\partial \chi} + v\frac{\partial u}{\partial \gamma} + w_{\sigma}\frac{\partial u}{\partial \sigma} = -\frac{1}{\rho} \left(\frac{\partial p_{D}}{\partial \chi} + \frac{\partial p_{D}}{\partial \sigma}\frac{\partial \sigma}{\partial x}  + \rho g\frac{\partial \eta}{\partial \chi}\right) + \nu \nabla_{\sigma}^2 u, \\
            & \frac{\partial v}{\partial \tau} + u\frac{\partial v}{\partial \chi} + v\frac{\partial v}{\partial \gamma} + w_{\sigma}\frac{\partial v}{\partial \sigma} = -\frac{1}{\rho} \left(\frac{\partial p_{D}}{\partial \gamma} + \frac{\partial p_{D}}{\partial \sigma}\frac{\partial \sigma}{\partial y}  + \rho g\frac{\partial \eta}{\partial \gamma}\right) + \nu \nabla_{\sigma}^2 v, \\
            & \frac{\partial w}{\partial \tau} + u\frac{\partial w}{\partial \chi} + v\frac{\partial w}{\partial \gamma} + w_{\sigma}\frac{\partial w}{\partial \sigma} = -\frac{1}{\rho} \frac{\partial p_{D}}{\partial \sigma}\frac{\partial \sigma}{\partial z} + \nu \nabla_{\sigma}^2 w. 
    \end{align}
\end{subequations}
Here the transformed vertical velocity $w_\sigma$ is introduced, which allows for the $\sigma$-transformed equations to be written in a form similar to the original equations.
The transformed gradient operator, the transformed Laplace operator, and the transformed vertical velocity are defined as follows
\begin{subequations}
\begin{gather}
        \nabla_{\sigma} = \left[\frac{\partial}{\partial \chi} + \frac{\partial \sigma}{\partial x}\frac{\partial}{\partial \sigma}, \frac{\partial}{\partial \gamma} + \frac{\partial \sigma}{\partial y}\frac{\partial}{\partial \sigma}, \frac{\partial \sigma}{\partial z}\frac{\partial}{\partial \sigma}\right]^T, \\
        \begin{aligned}
            & \nabla_{\sigma}^2 = \frac{\partial^2}{\partial \chi^2} + \frac{\partial^2}{\partial \gamma^2} + \left[ \left(\frac{\partial \sigma}{\partial x}\right)^2 + \left(\frac{\partial \sigma}{\partial y}\right)^2 + \left(\frac{\partial \sigma}{\partial z}\right)^2 \right]\frac{\partial^2 }{\partial \sigma^2} \\ &+ 2 \frac{\partial \sigma}{\partial x} \frac{\partial^2}{\partial \sigma \partial \chi}  + 2\frac{\partial \sigma}{\partial y} \frac{\partial^2}{\partial \sigma \partial \gamma}  + \left(\frac{\partial^2 \sigma}{\partial x^2} + \frac{\partial^2 \sigma}{\partial y^2} \right) \frac{\partial }{\partial \sigma},                
            \end{aligned}\\
            w_{\sigma} = \frac{\partial \sigma}{\partial t} + u\frac{\partial \sigma}{\partial x} + v\frac{\partial \sigma}{\partial y} + w\frac{\partial \sigma}{\partial z},
            \end{gather}
\end{subequations}
The $\sigma$-dependent derivatives are evaluated by applying the chain rule to the transformation function (\ref{eq:sigtrans}):
\begin{subequations} \label{eq:1_Scoef}
    \begin{gather}
        \frac{\partial \sigma}{\partial t} = -d^{-1}\left(\sigma \frac{\partial d}{\partial t}\right), \label{eq:1_ds/dt} \\
        \tilde{\nabla}\sigma = d^{-1}\left(\tilde{\nabla} h -\sigma \tilde{\nabla} d \right), \label{eq:1_ds/d}\\
        \tilde{\nabla}^2\sigma = d^{-1}\left(\tilde{\nabla}^2 h - \sigma \tilde{\nabla}^2 d - 2 \tilde{\nabla} \sigma \tilde{\nabla}d \right),  \\
        \frac{\partial \sigma}{\partial z} =  d^{-1}. 
    \end{gather} 
\end{subequations}
Note that $\tilde{\nabla}$ refers to differentiation in the horizontal plane $[\frac{\partial}{\partial x}, \frac{\partial}{\partial y}]^T$ only. To compute the time-dependent $\sigma$-coefficient, it is possible to use the kinematic boundary conditions
\begin{equation}
    \frac{\partial d\vert_{z=\eta}}{\partial t} = w\vert_{z=\eta} - u\vert_{z=\eta} \frac{\partial \eta}{\partial x}- v\vert_{z=\eta} \frac{\partial \eta}{\partial y}.
\end{equation}
Finally, the $\sigma$-transformed Navier-Stokes equations can be expressed as
\begin{subequations}
\begin{align}
    & \nabla_{\sigma} \cdot \textbf{u} = 0, \\
    &\frac{\partial \textbf{u}}{\partial t} + \textbf{u}_{\sigma}\cdot  \nabla \textbf{u} = - \frac{1}{\rho} (\nabla_{\sigma} p_D + \nabla p_S) + \textbf{g} + \nu \nabla_{\sigma}^2 \textbf{u},\label{eq:sigmom1} 
\end{align}
\end{subequations}
where $\nabla = \left[\frac{\partial}{\partial \chi},\frac{\partial}{\partial \gamma},\frac{\partial}{\partial \sigma}\right]^T$ and $\textbf{u}_\sigma =[u,v,w_\sigma]^T$.

\section{Numerical Discretization}
The general numerical discretization of the problem is done using a method of lines technique where the governing equations are discretized in space using spectral methods and then a standard ordinary differential equation (ODE) solver is applied on the resulting semi-discrete formulation for temporal integration.
\subsection{A General Low-Storage Runge-Kutta Method}
To advance the NSE in time, we use a $s$-stage low-storage explicit Runge-Kutta (LSERK) method as an ODE solver. Explicit low-storage Runge-Kutta methods have the general form \cite{Williamson80}
\begin{subequations}
\begin{align}
    &y^{(0)} = y^n,\quad K^0 = 0,\\
    &K^k = \alpha_kK^{k-1}+\Delta t f(y^{(k-1)},t_n+c_k\Delta t),\\
    &y^{(k)} = y^{(k-1)}+\beta_kK^k,\quad k=1,2,...,s, \label{eq:RKstage} \\
    &y^{n+1} = y^{(s)},
\end{align}
\label{eq:LSERK}
\end{subequations}
where $y^n$ refers to the solution at time $t_n$, $(k)$ refers to the stage, and $\alpha$, $\beta$ and $c$ are coefficients depending on a chosen LSERK method.

\subsection{Conservation of mass}\label{sec:conMass}
To ensure that the system fulfills conservation of mass, we have to ensure that $\nabla\cdot\textbf{u}^{(k)}=0$ at all RK stages. To do this a pressure-correction method is used, where a Poisson boundary value problem is solved at each stage in the Runge-Kutta method. The momentum equations (\ref{eq:sigmom1}) is denoted as
\begin{align}
   f(\textbf{q}) =-\textbf{u}_\sigma\cdot\nabla \textbf{u}-\frac{1}{\rho}(\nabla_\sigma p_D +\nabla p_S)  + g + \nu \nabla_\sigma^2 \textbf{u}\label{eq:forc1},
\end{align}
with $\textbf{q}=[\textbf{u},p,\eta]$. Combining this with the general LSERK method (\ref{eq:LSERK}) we get
\begin{align}
    \textbf{u}^{(k)} &= \textbf{u}^{(k-1)}+\beta_{k}K^k\nonumber\\
    &= \textbf{u}^{(k-1)}+\beta_{k}(\alpha_kK^{k-1}+\Delta t f(\textbf{q}^{(k-1)})).\label{eq:rkforc1}
\end{align}
To find a discrete pressure-velocity coupling similar to \eqref{eq:velocitypressurecont}, 
the divergence at a new RK stage $k$ can be found by applying the gradient operator to \eqref{eq:RKstage} 
\begin{align}\nabla_\sigma^k\cdot\textbf{u}^{(k)} =\nabla_\sigma^k\cdot\textbf{u}^{(k-1)}+\nabla_\sigma^k\cdot\beta_{k}(\alpha_kK^{k-1}+\Delta t f(\textbf{q}^{(k-1)})).\label{eq:rkdiv1}
\end{align}
Writing out the last forcing term gives the following
\begin{align}
    &\nabla_\sigma^k\cdot\textbf{u}^{(k)} = \nabla_\sigma^k\cdot\textbf{u}^{(k-1)}+\nabla_\sigma^k\cdot\beta_{k}\alpha_kK^{k-1}\nonumber \\
    &+ \nabla_\sigma^k\cdot\beta_{k}\Delta t (-\frac{1}{\rho}(\nabla_\sigma^{k-1} p_D^{(k-1)} +\nabla^{k-1} p_S^{(k-1)}) \nonumber\\& + g + \nu (\nabla_\sigma^{k-1})^2 \textbf{u}^{(k-1)}-\textbf{u}_\sigma^{(k-1)}\cdot\nabla\textbf{u}^{(k-1)}).
\end{align}
It is now possible to isolate the dynamic pressure, $p_D$, on the left hand side, and utilize that $\nabla^{k}_{\sigma}\cdot\textbf{u}^{(k)}=0$, leading to a mixed-stage Poisson pressure problem, defined as
\begin{align}&\nabla_\sigma^k\cdot\nabla_\sigma^{k-1}p_D^{(k-1)} = \frac{\rho}{\beta_k\Delta t}\nabla_\sigma^k\cdot\textbf{u}^{(k-1)}+\frac{\rho\alpha_k}{\Delta t}\nabla_\sigma^k\cdot K^{k-1}\nonumber\\
    &-\nabla_\sigma^k\cdot\nabla^{k-1} p_S^{(k-1)}  + \rho\nabla_\sigma^k\cdot g + \rho\nu\nabla_\sigma^k\cdot (\nabla_\sigma^{k-1})^2 \textbf{u}^{(k-1)}\nonumber \\
    &-\rho\nabla_\sigma^{k}\cdot(\textbf{u}_\sigma^{(k-1)}\cdot\nabla\textbf{u}^{(k-1)}).\label{eq:poissonp}
\end{align}
We now have derived a mixed-stage Poisson problem for the dynamic pressure $p_D$ that ensures a divergence free velocity field, however, the boundary conditions still have to be defined for the problem to be well-posed.

\subsection{Boundary conditions for pressure}\label{sec:bc}
As mentioned earlier the boundary condition for the pressure at the surface is given as
\begin{align}
 p_D = 0 \quad \text{on}\quad z=\eta(x,y,t)\label{eq:bcpd}.
\end{align}
The boundary condition of zero pressure at the free surface is a Dirichlet boundary condition. 
At solid boundaries, impermeability boundary conditions are applied, given as
\begin{align}
\textbf{n}\cdot\textbf{u} = 0\quad \text{on}\quad\Gamma_w.\label{eq:bcun}   
\end{align}
The impermeable condition for the velocity can be transformed to a dynamic pressure condition by taking the normal to the LSERK expression in (\ref{eq:rkforc1}). The derivation then follows the same way as for the conservation of mass. This gives the impermeability boundary condition expressed as a boundary condition for the dynamic pressure Poisson problem,
\begin{align}
    &\textbf{n}\cdot\nabla_\sigma^{k-1}p_D^{(k-1)} = \frac{\rho}{\beta_k\Delta t}\textbf{n}\cdot\textbf{u}^{(k-1)}+\frac{\rho\alpha_k}{\Delta t}\textbf{n}\cdot K^{k-1}\nonumber\\
    &-\textbf{n}\cdot\nabla^{k-1} p_S^{(k-1)}  + \rho\textbf{n}\cdot g + \rho\nu\textbf{n}\cdot (\nabla_\sigma^{k-1})^2 \textbf{u}^{(k-1)}-\rho\textbf{n}\cdot(\textbf{u}_\sigma^{(k-1)}\cdot\nabla\textbf{u}^{(k-1)}).\label{eq:bcpdneu}
\end{align}
Now, \eqref{eq:poissonp} in combination with \eqref{eq:bcpd} and \eqref{eq:bcpdneu} defines a well-posed Poisson problem that can be solved to determine the dynamic pressure $p_D$ and a divergence free velocity field can be obtained through updating the evolution equations \eqref{eq:kinematicFS} and  \eqref{eq:sigmom1} using the general LSERK scheme \eqref{eq:LSERK}. 

\subsection{Spatial discretization in 2D}\label{sec:spat1}
We discretize the NSE using a spectral collocation method. The equations are represented by a mixed basis consisting of a Fourier collocation method in the horizontal $x$-direction, and a Chebyshev collocation method in the vertical $\sigma$-direction.
The Fourier discretization is based on the Fourier basis  $\exp(ikx)$, with $k = -\frac{N}{2},...,\frac{N}{2}$ being the Fourier wave numbers and $N$ being the order of the discretization. As the notation for even and uneven $N$ differs slightly, we will here assume that $N$ is even to simplify the notation, but note that the method is valid for both even and uneven. The nodes for the horizontal grid are distributed equidistantly as 
\begin{align}
    x_n = \frac{2\pi n}{N+1},\quad n=0,..,N.
\end{align}
The Chebyshev discretization in the vertical $\sigma$-direction, is based on the Chebyshev polynomials which are given as
\begin{align}
    T_m(\sigma) = \cos(m\arccos(\sigma)),\quad m\geq0,\quad -1\leq\sigma\leq1.
\end{align}
The nodes for the horizontal grid are distributed according to the Chebyshev Gauss-Lobatto point distribution as
\begin{align}
    \sigma_m =\cos\left(\frac{m\pi}{M+1}\right),\quad m = 0,..,M,
\end{align}
with $M$ being the order of the discretization.
Based on these we can write the discretization of the free surface, the velocities, and the dynamic pressure as
\begin{subequations}
\begin{gather}
    \begin{aligned}
    \eta(x,\sigma,t) &=  \sum_{n=-\frac{N}{2}}^{\frac{N}{2}} \hat{\eta}_{n,m}(t)\exp(ik_nx)\\&=\sum_{n=-\frac{N}{2}}^{\frac{N}{2}} \eta_{n,m}(x_n,\sigma_m,t)l_n(x),
    \end{aligned}\\
    \begin{aligned}
    \textbf{u}(x,\sigma,t) &=  \sum_{n=-\frac{N}{2}}^{\frac{N}{2}}\sum_{m=0}^{M} \hat{\textbf{u}}_{n,m}(t)\exp(ik_nx)T_m(\sigma)\\&=\sum_{n=-\frac{N}{2}}^{\frac{N}{2}}\sum_{m=0}^{M} \textbf{u}_{n,m}(x_n,\sigma_m,t)l_n(x)l_m(\sigma),
    \end{aligned}\\
    \begin{aligned}
    p_D(x,\sigma,t) &=  \sum_{n=-\frac{N}{2}}^{\frac{N}{2}}\sum_{m=0}^{M} \hat{p}_{Dn,m}(t)\exp(ik_nx)T_m(\sigma)\\&=\sum_{n=-\frac{N}{2}}^{\frac{N}{2}}\sum_{m=0}^{M} p_{Dn,m}(x_n,\sigma_m,t)l_n(x)l_m(\sigma).
    \end{aligned}
\end{gather}
\end{subequations}
Here "\^{}" denotes the modal trigonometric coefficients, while no hat denotes the nodal coefficients. Moreover, $l_n(x)$ and $l_m(\sigma)$ are the nodal Lagrange interpolation polynomials for the Fourier and Chebyshev representation respectively, given as $l_n(x_m) = \delta_{n,m}$, where $\delta_{n,m}$ is Kroneckers delta defined as
\begin{align}
\delta_{n,m} = \left\{ \begin{array}{cc} 
                0 & \hspace{5mm} n\neq m,\\
                1 & \hspace{5mm} n = m. \\
                \end{array} \right.
\end{align} 
\subsection{Computation of derivatives}\label{sec:deriv}
The Fourier and Chebyshev discretization benefit from being able to compute derivatives in $O(Nlog(N))$ time through the use of the Fast Fourier transform (FFT), $\mathcal{F}$, and the Fast Chebyshev transform (FCT), $\mathcal{F}_c$, respectively. 
For the Fourier discretization, we can evaluate the continuous derivatives $\partial_x$ and $\partial_{xx}$ by employing the FFT, to find the Fourier coefficients. It is then well known, that the first and second derivative can be found by multiplying the coefficients by $ik$ and $-k^2$ for the first and second derivative respectively. Afterwards the inverse Fast Fourier Transform (IFFT), $\mathcal{F}^{-1}$, can be employed to return to the nodal space. By letting $\hat{D}_1$ and $\hat{D}_2$ denote the multiplication by $ik$ and $-k^2$ respectively, the computation of the horizontal derivatives can be written as 
\begin{subequations}
\begin{align}
    \partial_xf &\approx D_1f = \mathcal{F}^{-1}(\hat{D}_1\mathcal{F}(f)),\\
     \partial_{xx}f &\approx D_2f = \mathcal{F}^{-1}(\hat{D}_2\mathcal{F}(f),
\end{align}
\end{subequations}
and the same operations apply to the $y$-direction.

For the Chebyshev discretization in the vertical direction, the derivatives can be evaluated by either evaluating the derivatives of the Chebyshev polynomials or their modal coefficients as
\begin{align}
\label{eq:chebDeriv}
    \frac{\partial^{(p)}}{\partial\sigma^{(p)}}f_{n,m} = \sum_{m=p}^{M}\hat{f}_{n,m}\frac{\partial^{(p)}}{\partial\sigma^{(p)}}T_m(\sigma) = \sum_{m=0}^{M-p}\hat{f}_{n,m}^{(p)}T_m(\sigma),
\end{align}
with $\hat{f}_{n,m}$ being the modal Chebyshev coefficients.
The modal coefficients related to the basis of the solutions or derivative hereof are related through the following recursive formula
\begin{align}
    \hat{f}^{(p)}_{n,q-1} = \frac{1}{c_{q-1}}\left(\hat{f}^{(p)}_{n,q+1}+2q\hat{f}^{(p-1)}_{n,q}\right).
\end{align}
The recursion is initialised with $\hat{f}^{(p)}_{n,M-p+1}=0$, with $q$ running backwards as $q = M-p,...,1$ with $c_q = 1+\delta_{q,0}$. Through the recursion the first and second derivative of the modal coefficients, $\hat{f}^{(1)}_{n,m}$ and $\hat{f}^{(2)}_{n,m}$, can be evaluated in $O(M)$ time. In combination with utilizing FCT, $\mathcal{F}_c$, to evaluate the modal Chebyshev coefficients along with (\ref{eq:chebDeriv}), the nodal derivatives can then be evaluated in $O(M\log M)$ time. For any mixed derivatives, e.g. $\partial_{x\sigma}$, the two methods are simply applied one after the other.

\subsection{Dealing with aliasing stemming from nonlinear terms}
Nonlinear terms occur in the governing equations, i.e. the momentum equations \eqref{eq:moms2}, the kinematic free surface boundary condition \eqref{eq:kinematicFS}, and the right hand side terms of the pressure Poisson problem \eqref{eq:poissonp}, \eqref{eq:bcpd} and \eqref{eq:bcpdneu}. It is well known that such terms can introduce aliasing errors when evaluated, which in term can affect accuracy and stability of the model. Aliasing errors occur due to the interpolated product of functions having higher order than the used basis. This leads to high-frequency errors polluting the higher-order modes of the solution if untreated. To mitigate effects of aliasing-driven errors that may lead to numerical  instabilities, we introduce an  exponential cutoff filter defined as
\begin{equation}
     F_S(i,N_c,\alpha,s) = \begin{cases} 
          1 & 0\leq i\leq N_c \\
          \exp\left(\alpha\left(\frac{i-N_c}{N+1-N_c}\right)^s\right) & N_c<i\leq N 
       \end{cases}.
\end{equation}
The filter exponentially reduces the energy of modes higher than the cut-off order $N_c$, with $\alpha$ and $s$ being parameters. The filter is applied independently to the Fourier and Chebyshev expansions. To define a mild low-pass filter designed to not destroy the high-order accuracy of the numerical scheme, we use $N_c = 0.5N_x$ and $N_c=0.9N_\sigma$, respectively, with parameters $\alpha=36$ and $s=2$. The filter is applied on the modes at every stage of the LSERK through FFT and FCT, 
\begin{subequations}
\begin{align}
    f^{F}_{filtered} &= \mathcal{F}^{-1}(F_S(\mathcal{F}(f))),\\
    f^{C}_{filtered} &= \mathcal{F}_c^{-1}(F_S(\mathcal{F}_c(f))),
\end{align}
\end{subequations}
for the Fourier and Chebyshev expansions, respectively. 

\section{Geometric $p$-Multigrid Method}
\begin{figure}
    \centering
    \includegraphics[width=70mm]{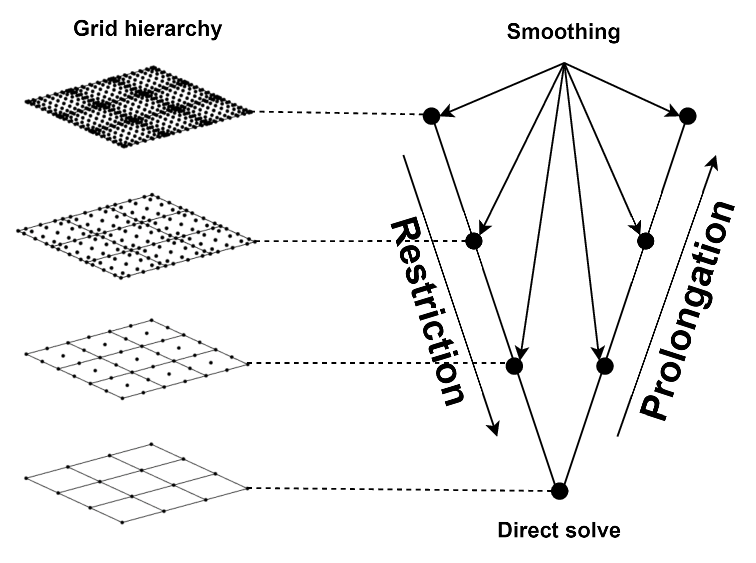}
    \caption{Example of a Geometric multigrid V-cycle where the grid hierarchy is traversed through utilization of transfer operators.}
    \label{fig:vcycle}
\end{figure}

Designing an efficient iterative method to solve the occurring Poisson pressure problem is an important part of solving the NSE with a free surface. We are interested in utilizing iterative methods for two main reasons: Iterative methods outperform direct solvers for large systems, and they let us avoid having to explicitly form the system matrix $\boldsymbol{A}$, so we can take advantage of FFT and FCT as mentioned in section \ref{sec:deriv}. Among the iterative methods, the multigrid class is often the most efficient linear system of equations solver and the $p$-multigrid methods are less explored. In this part we employ a geometric $p$-multigrid method to solve the matrix system resulting from the high-order discretization models. The $p$-multigrid method takes advantage of the high-order spatial discretization, as this allows the coarsening of grids to be done by simply reducing the polynomial order $p$. Moreover, the modal representation of the solutions allows for easy building of transfer operators necessary to move the solutions between meshes while maintaining the numerical advantages of the spectral accurate discretization. Figure \ref{fig:vcycle} shows an example of a V-cycle $p$-multigrid method. The V-cycle is the simplest multigrid cycle, and can be seen as the "building block" of any other types of cycles. The solution starts on the fine grid, where smoothing is applied. The solution is then moved to a coarser grid using a restriction operator, before smoothing is applied again. This continues until the solution reaches the coarsest grid. Here a direct solver is often used instead of simply smoothing, where afterwards the solution is moved up through the fine grids again using a prolongation operator, with smoothing again being applied at every grid level. 
In this paper we will focus on the $V$-cycle and the Full Multigrid Method (FMG), which can be seen as a successive application of $V$-cycles. The pseudo algorithm for the $V$-cycle and FMG can be seen in Algorithm \ref{alg:Vcycle} and Algorithm \ref{alg:FMG} respectively.

\begin{algorithm}[H]
\caption{Recursive Multigrid $V$-cycle algorithm}\label{alg:Vcycle}
\begin{algorithmic}[1]
\Function{Vcycle}{$p,A,u,f,v_1,v_2$}
\If{$p=1$}
    \State Solve system using a direct solver on coarsest grid level: $u=A^{-1}f$
    \State return $u$
\EndIf
\State Apply pre-smoothing $v_1$ times
\State Compute residual: $r = f-Au$
\State Restrict residual to coarser grid: $r_c = R(r)$
\State Set initial state of error to zero: $e_c = 0$
\State Recursive call to $V$-cycle: Vcycle($p-1,A,e_c,r_c,v_1,v_2$)
\State Prolongate correction to fine grid: $e=P(e_c)$
\State Correct solution: $u=u+e$
\State Apply post-smoothing $v_2$ times
\State return $u$
\EndFunction
\end{algorithmic}
\end{algorithm}

\begin{algorithm}[H]
\caption{FMG}\label{alg:FMG}
\begin{algorithmic}[1]
\Function{FMG}{$A,u,f,v_1,v_2$}
\For{p=1:FinestGrid-1}
\State Call $V$-cycle on current grid level: $u=\text{Vcycle}(p,A,u,f,v_1,v_2)$
\State Prolongate result to next grid level: $u=P(u)$
\EndFor
\State Call $V$-cycle on finest grid level: $u=\text{Vcycle}(p,A,u,f,v_1,v_2)$
\State return $u$
\EndFunction
\end{algorithmic}
\end{algorithm}

\subsection{Grid coarsening strategy}
The coarsening of grids is done by establishing a grid hierarchy obtained by  reducing the polynomial order $p$ of the spatial discretization. We let $p=1,..,P$ be the levels of the multigrid method, and consider the grids $G^p$ on the domain $\Omega$ with order $N^p$ and $M^p$ in the horizontal and vertical direction respectively, with $G^P$ being the finest grid, and $G^1$ the coarsest. We let the order of discretization of the grids for the combined standard- and semi-coarsening strategy be 
\begin{align*}
  N^p &=\Bigg\{\begin{array}{ll}
      2N^{p-1} & \texttt{for} \quad N^{p-1}=M^{p-1} \texttt{ and } N^{p-1}<M^{p-1}\\
      N^{p-1}  & \texttt{otherwise}
    \end{array},\\
  M^p &=\Bigg\{\begin{array}{ll}
      2M^{p-1} & \texttt{for} \quad M^{p-1}=N^{p-1} \texttt{ and } M^{p-1}<N^{p-1}\\
      M^{p-1} & \texttt{otherwise}
    \end{array}.\\
\end{align*}
This gives us a general method for constructing a group of nested grids. While there are no strict requirements for the resolution of grid $G^1$, it is generally chosen to be "coarse enough" such that a direct solve can be done in reasonable time.

\subsection{Transfer operators}
An important aspect of the geometric multigrid method is being able to efficiently move solutions between grids. Moreover, we have to ensure that we conserve the spectral convergence properties while transferring. This is however rather simple due to the modal representations of our solutions introduced in section \ref{sec:spat1}. We have to define both a restriction operator, which can move the solution from fine to coarse grid $R^p: G^{p-1} \rightarrow G^p$ and a prolongation operator that can move the solution from coarse to fine grid $P^p: G^{p} \rightarrow G^{p-1}$. The transfer in the horizontal is done using FFT. By using FFT to compute the Fourier modes of the solution, it is simple then to transfer from one grid to another. When moving from fine to coarse, i.e. restriction, we simply remove higher order modes before using IFFT to move back to nodal values. For prolongation, we again use FFT to compute the Fourier modes, and then zero-pad before using IFFT to move back to nodal values. The same course of action can be done with the vertical discretization. We simply employ FCT to compute the modal coefficient instead, and otherwise follow the same method. 

\subsection{Smoothing strategy}
Smoothing consists of applying an iterative solver to the problem at any grid level, with the solver typically being in the family of stationary iterative solvers. This iterative solver for smoothing is not applied until convergence, as this would take away the benefits of the multigrid method, but is instead only applied for some amount of iterations. This is what is responsible for the actual error correction, and can have large impact on the iterations required for convergence. Here we employ a blocked Gauss-Seidel method, where the computational domain is divided across the horizontal dimensions into $K$ blocks. The blocked Gauss-Seidel method generally has faster convergence, and avoids explicitly constructing the full system matrix, but instead solves a series of smaller systems.
\begin{align}
    \textbf{x}_r^{(m+1)} = \textbf{A}_{r,r}^{-1}\left(\textbf{f}_r-\sum_{s=1}^{r-1}\textbf{A}_{r,s}\textbf{x}_s^{(m+1)}-\sum_{s=r+1}^{N_\sigma}\textbf{A}_{r,s}\textbf{x}_s^{(m)}\right), \quad m=1,..,K-1,
\end{align}
where $\textbf{A}_{i,i}$ are $(N_\sigma\times N_\sigma)$ square matrices corresponding to the system matrices for the block of values $\textbf{x}_r$. 

\subsection{Combination with a Krylov-subspace method}
While it is possible to solve the full Poisson problem solely using a geometric $p$-multigrid method, it is often more robust to combine it with another iterative method. In this paper we therefore combine the $p$-multigrid method with a Krylov subspace method such as the GMRES method \cite{Saad86}, by utilizing the multigrid method as a preconditioner similar to the accelerated solver proposed for FNPF modeling using the spectral element method in \cite{AllanLaskowski}. This is done by applying the $p$-multigrid for a single cycle in the preconditioning step of GMRES. This has the benefit of reducing the condition number of the system matrix, and thereby the iteration count and overall computational cost. 

\section{Results}
\subsection{Geometric $p$-multigrid efficiency}
Table \ref{tab:PMGtime} shows the computation times for the $V$-cycle $p$-multigrid preconditioned GMRES method, the FMG $p$-multigrid preconditioned GMRES method and the GMRES method with no preconditioning with spatial polynomial orders in a 2D setup where polynomial orders are $N_x=40$ and $N_\sigma=40$. The solved problem is the Poisson problem derived in section \ref{sec:conMass} and considered for a stream function wave \cite{Rienecker_Fenton_1981} with periodic boundary conditions in the horizontal, and surface and bottom boundary conditions as defined in section \ref{sec:bc}. As expected the geometric $p$-multigrid preconditioner clearly reduces both iteration count and computational time for the Poisson solve. While the FMG leads to slightly fewer iterations than the V-cycle, the V-cycle still wins out in computational time, due to the lesser cost per iteration compared to FMG.
\begin{table}[H]
\begin{center}
{\scriptsize
\resizebox{\columnwidth}{!}{%
\begin{tabular}{cccccccccccc}
\hline
\multicolumn{3}{|c|}{}&
 \multicolumn{3}{|c|}{GMRES}&
 \multicolumn{3}{c|}{V-Cycle+GMRES}&
 \multicolumn{3}{c|}{FMG+GMRES}\\
\hline
\multicolumn{3}{|c|}{Rel. Error:}&
 \multicolumn{1}{|c|}{$10^{-4}$}&
 \multicolumn{1}{c|}{$10^{-8}$}&
 \multicolumn{1}{c|}{$10^{-12}$}&
 \multicolumn{1}{|c|}{$10^{-4}$}&
 \multicolumn{1}{c|}{$10^{-8}$}&
 \multicolumn{1}{c|}{$10^{-12}$}&
 \multicolumn{1}{|c|}{$10^{-4}$}&
 \multicolumn{1}{c|}{$10^{-8}$}&
 \multicolumn{1}{c|}{$10^{-12}$}\\
 \hline
\multicolumn{3}{|c|}{Iterations:}&
 \multicolumn{1}{|c|}{542}&
 \multicolumn{1}{c|}{987}&
 \multicolumn{1}{c|}{1183}&
 \multicolumn{1}{c|}{2}&
 \multicolumn{1}{c|}{4}&
 \multicolumn{1}{c|}{5.5}&
 \multicolumn{1}{c|}{1.75}&
 \multicolumn{1}{c|}{3.25}&
 \multicolumn{1}{c|}{5}\\
 \hline
 \multicolumn{3}{|c|}{Time [s]:}&
 \multicolumn{1}{|c|}{0.81}&
 \multicolumn{1}{c|}{2.25}&
 \multicolumn{1}{c|}{3.16}&
 \multicolumn{1}{c|}{0.26}&
 \multicolumn{1}{c|}{0.42}&
 \multicolumn{1}{c|}{0.58}&
 \multicolumn{1}{c|}{0.34}&
 \multicolumn{1}{c|}{0.51}&
 \multicolumn{1}{c|}{0.66}\\
 \hline
\end{tabular}
}
}%
\caption{Computational times for solving the pressure Poisson problem \eqref{eq:poissonp}.}
\label{tab:PMGtime}
\end{center}
\end{table}

\subsection{Nonlinear accuracy}
\begin{figure}[H]
\vspace*{0cm}
\hspace*{0cm}
\centering
\subfloat[]{
  \includegraphics[width=38mm]{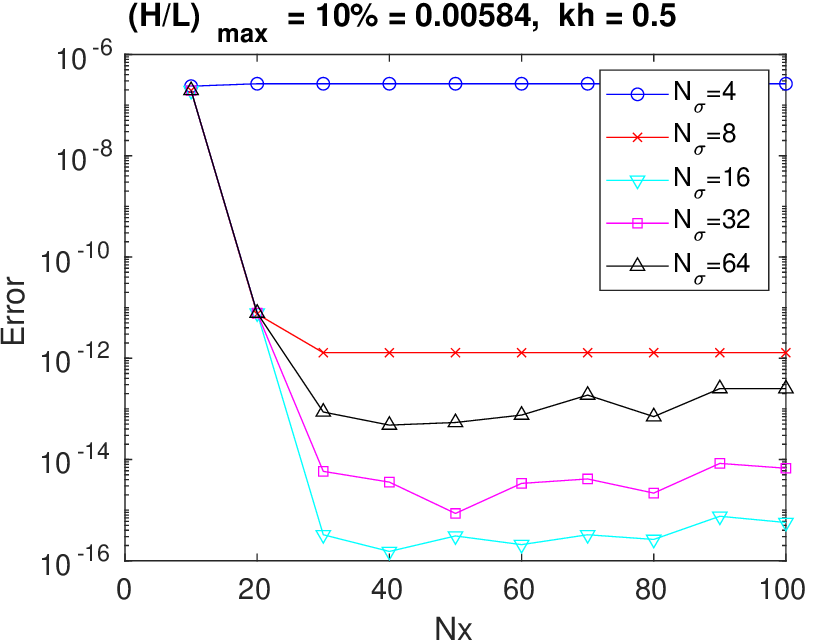}
}
\subfloat[]{
  \includegraphics[width=38mm]{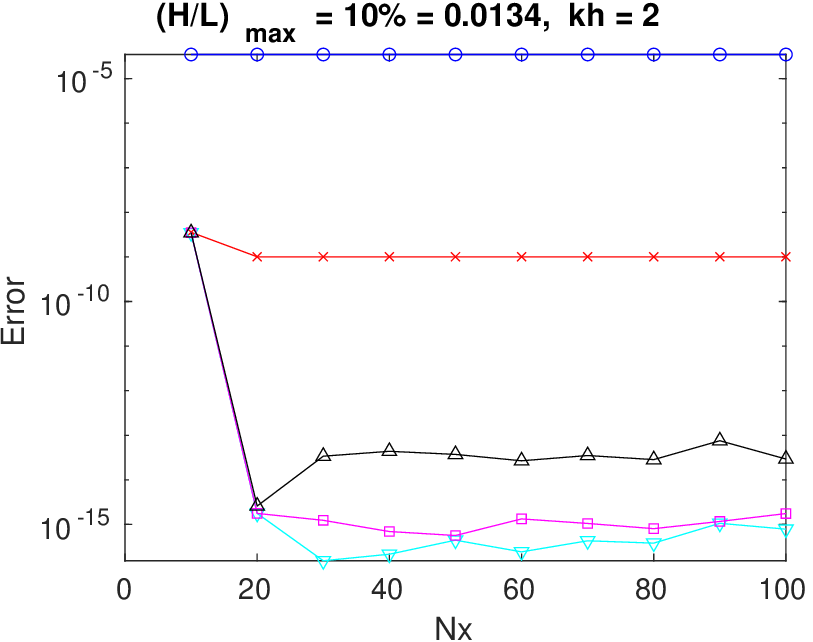}
}
\subfloat[]{
  \includegraphics[width=38mm]{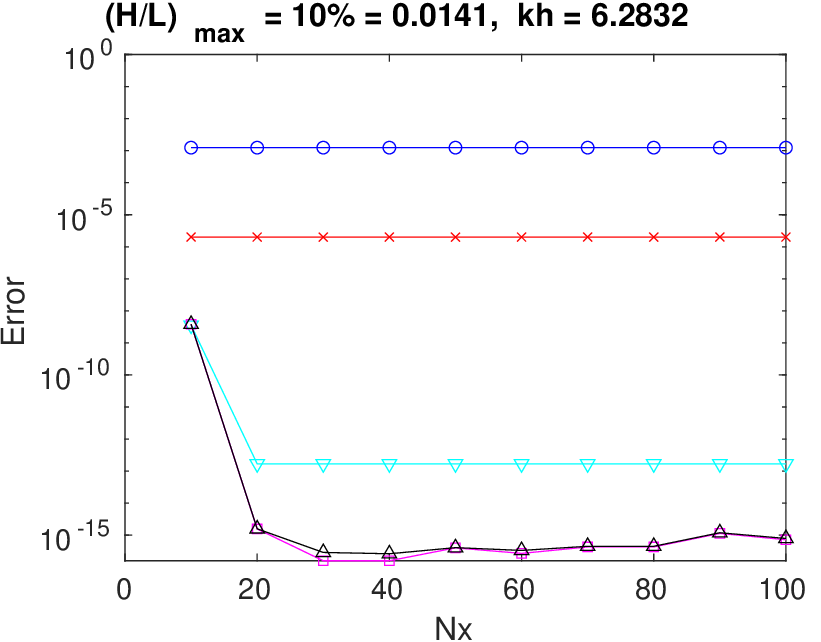}
}
\hspace{0mm}
\hspace*{0cm}
\subfloat[]{
  \includegraphics[width=38mm]{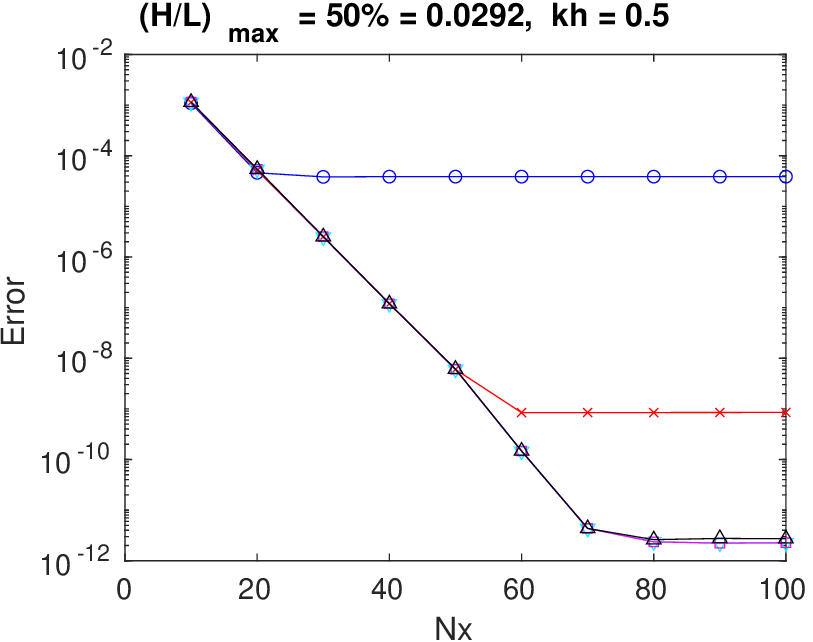}
}
\subfloat[]{
  \includegraphics[width=38mm]{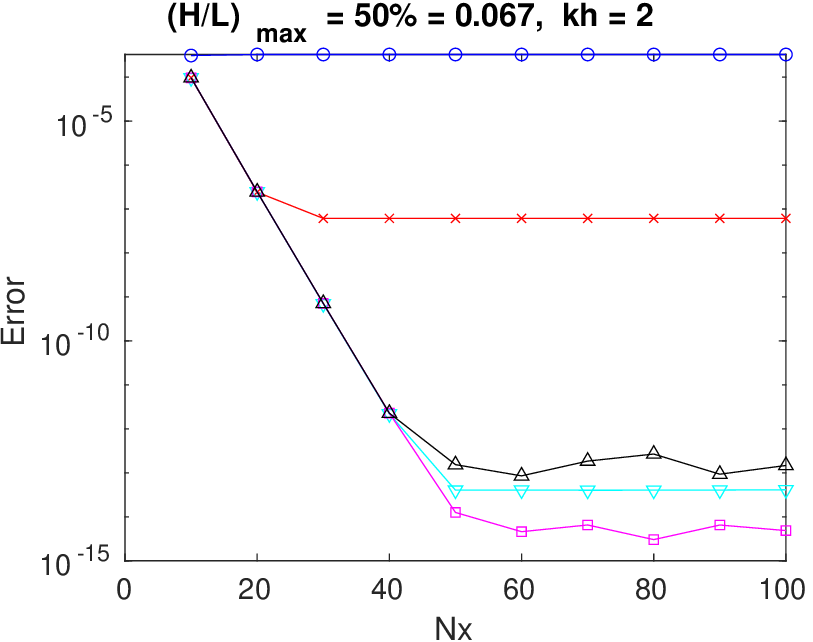}
}
\subfloat[]{
  \includegraphics[width=38mm]{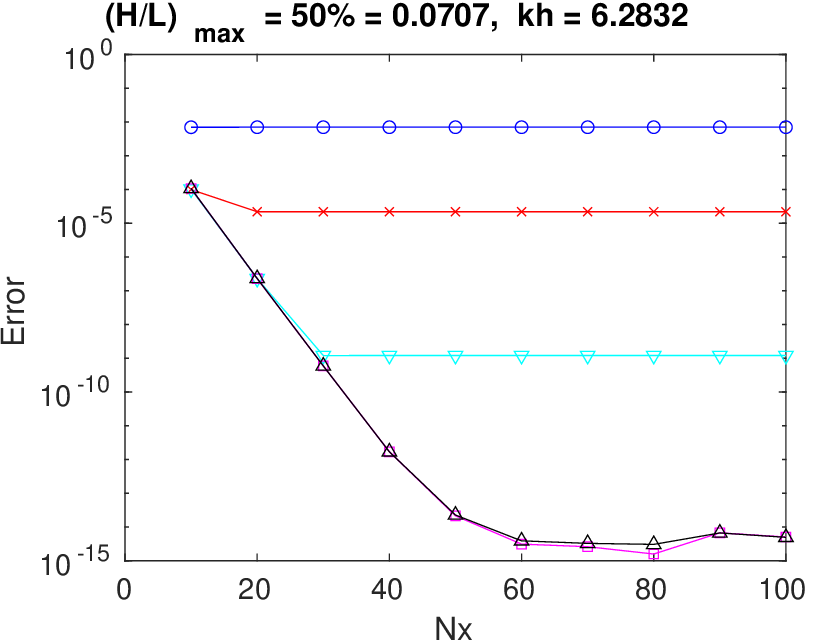}
}
\hspace{0mm}
\hspace*{0cm}
\subfloat[]{
  \includegraphics[width=38mm]{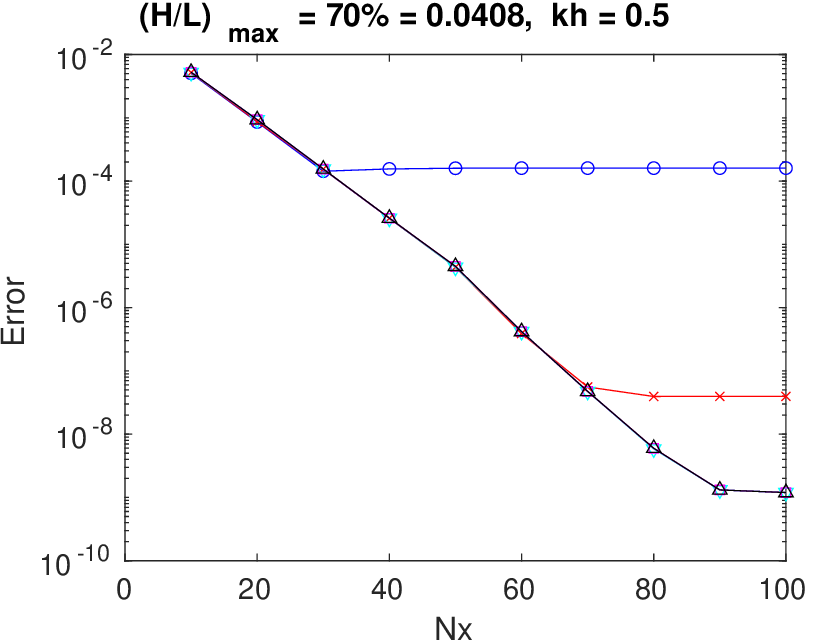}
}
\subfloat[]{
  \includegraphics[width=38mm]{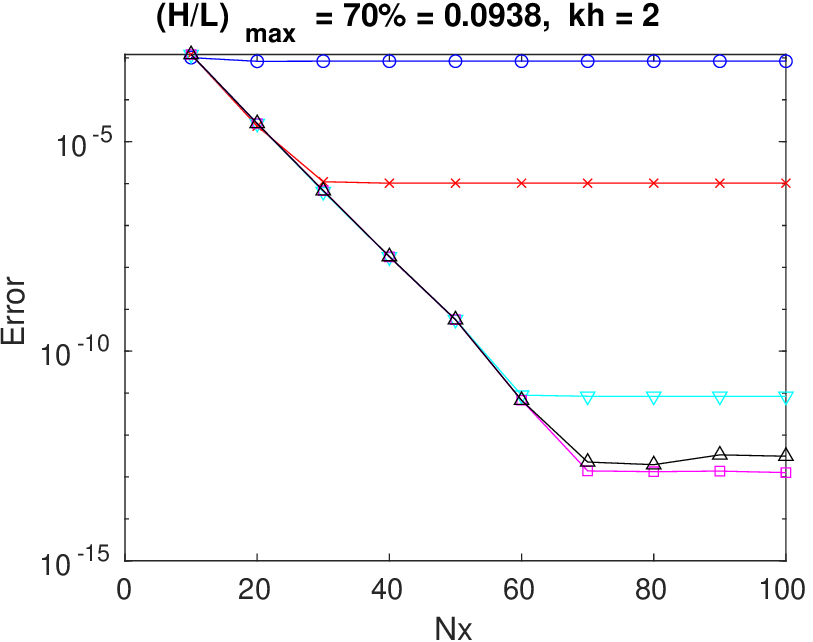}
}
\subfloat[]{
  \includegraphics[width=38mm]{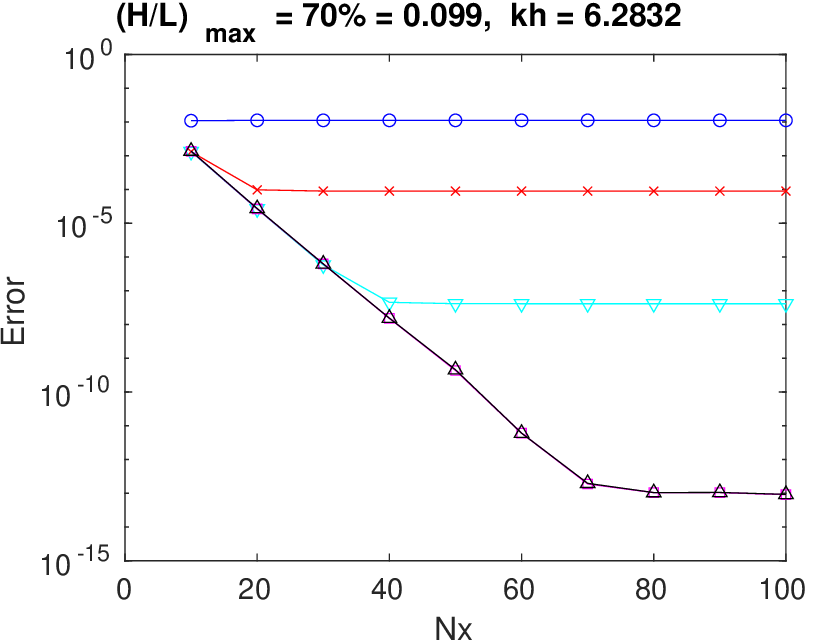}
}
\hspace{0mm}
\hspace*{0cm}
\subfloat[]{
  \includegraphics[width=38mm]{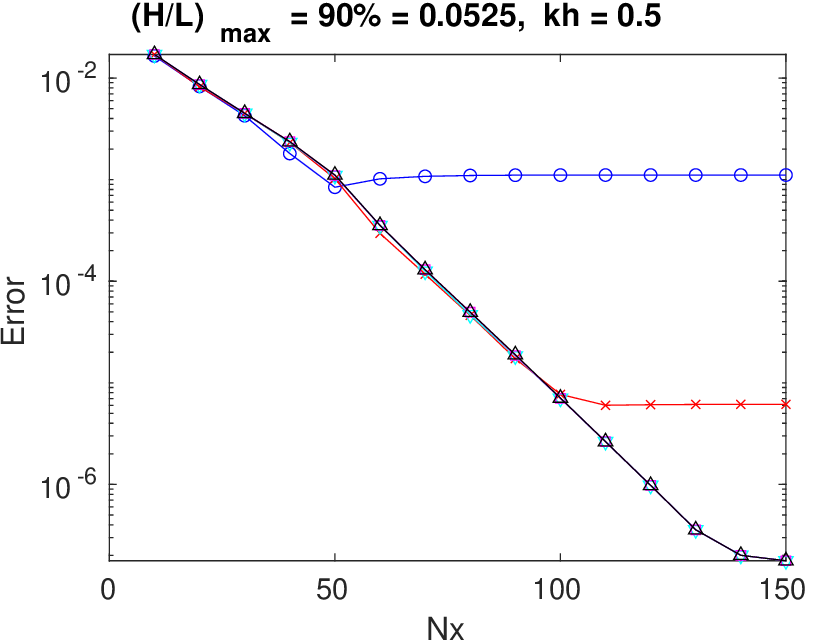}
}
\subfloat[]{
  \includegraphics[width=38mm]{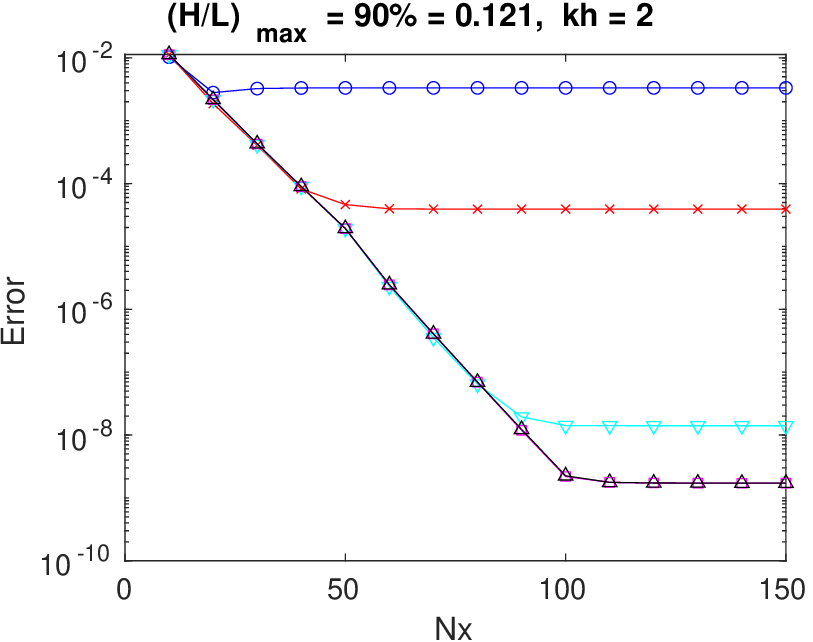}
}
\subfloat[]{
  \includegraphics[width=38mm]{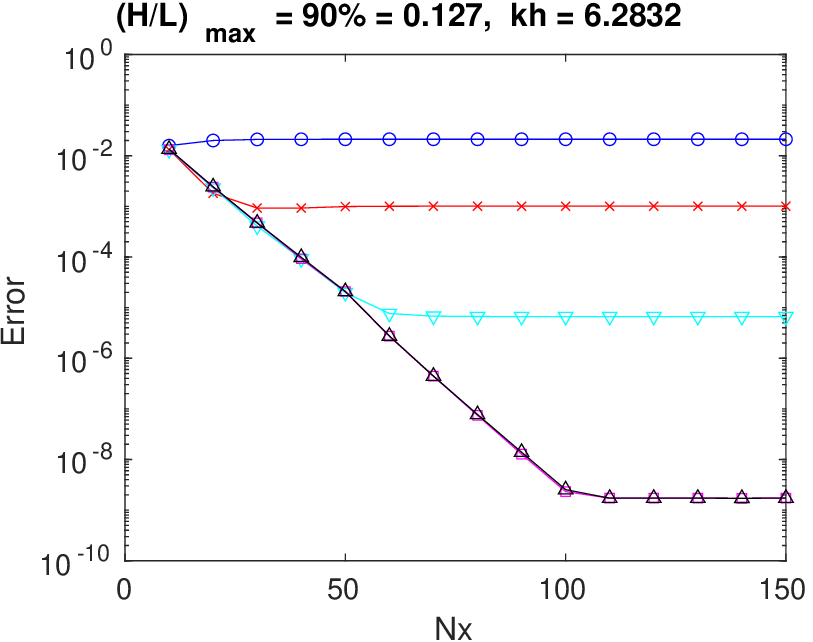}
}
\caption{Convergence plots for $u$ for varying nonlinearity ($H/L$) across rows in the columns and varying dispersion ($kh$) across columns in each row.}
\label{fig:uconv}
\end{figure}

To highlight the attractive spectral accuracy of the spatial discretization and to validate our model, we present a convergence study of the velocities. The waves are initialized by the stream function solution, which gives an analytical solution to inviscid, irrotational flow. The study is computed by computing the error after a single time step, such that the error of both the Poisson equation and the momentum equations are taken into account. The time step is set to be small enough such that spatial errors dominate. The study has been done on shallow ($kh=0.5$), intermediate ($kh=2$), and deep water ($kh=2\pi$) for four cases of non-linear waves. Here $kh$ denotes the non-dimensional wavenumber comprised of the still water depth $h$ and the wavenumber $k$. The wave steepness is defined as the maximum allowed steepness before breaking, as determined by Battjes \cite{Battjes74} 
\begin{align}
    \left(\frac{H}{L}\right)_{max} = 0.1401 \tanh(0.8863kh),
\end{align}
with $H$ being the wave height [m], and $L$ the wave length [m].
The resulting twelve convergence tests can be seen in Figure \ref{fig:uconv}. 

As can be seen, the method shows the expected spectral convergence for all cases. We note that shallow water in general requires higher resolution in the horizontal, while deep water requires increased vertical resolution to obtain high accuracy. Moreover, increased wave non-linearity ($H/L$) also requires higher resolution in both the horizontal and vertical directions to obtain high accuracy solutions.

\subsection{Linear dispersion analysis}
\begin{figure}[H]
\vspace*{0cm}
\hspace*{0cm}
\centering
\includegraphics[width=80mm]{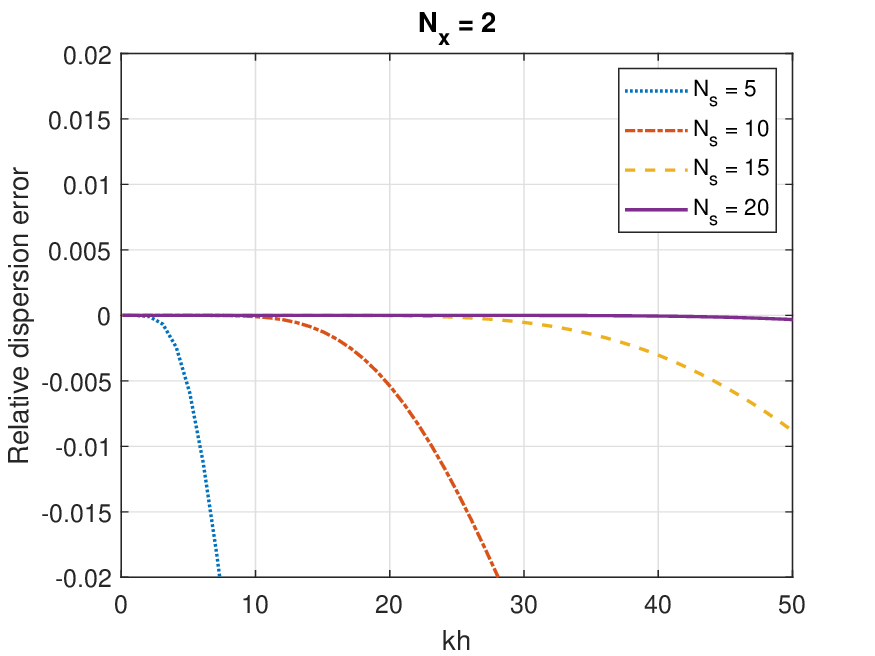}
\caption{Linear dispersion error for $N_x = 2$.}
\label{fig:LDE}
\end{figure}

To show the models ability to accurately simulate waves and flow kinematics with low dispersion, we present a linear dispersion analysis. For a periodic linear solution, the vertical surface velocity $\tilde{w}$ and potential $\tilde{\phi}$ are related through the linear dispersion operator \cite{Allan09}
\begin{align}
    \frac{\tilde{w}}{k\tilde{\phi}} = \tanh(kh),
\end{align}
which means the simulated potential can be computed from the vertical free surface velocity as
\begin{align}
    \tilde{\phi_s}=\frac{\tilde{w}_s}{k\tanh(kh)}.
\end{align}
From this, the relative dispersion error (RDE) can be computed by comparing the exact and simulated potential
\begin{align}
    RDE = \frac{\|\phi_e-\phi_s\|_2}{\|\phi_e\|_2},
\end{align}
with the exact solution from potential flow theory as
\begin{align}
    \phi_e = -\frac{Hc}{2}\frac{\cosh(k(z+h))}{\sinh(kh)}
\sin(\omega t-kx).
\end{align}
To carry out this analysis, we consider the linearised Navier-Stokes equations under the assumption that $\frac{H}{L}\ll 1$, i.e. small amplitude. Under this assumption, the governing equations are reduced to
\begin{subequations}
    \begin{align}
        & \nabla \cdot \textbf{u} = 0, \\ 
        &\frac{\partial \textbf{u}}{\partial t} = - \frac{1}{\rho} \nabla p + \textbf{g} +\nu\nabla^2\textbf{u},\label{eq:momslin}
    \end{align} \label{eq:NSElin}%
\end{subequations}
with the kinematic free surface boundary condition reduced to
\begin{equation}
    \frac{\partial\eta}{\partial t} =\tilde{w}.
    \label{eq:linFS}
\end{equation}
Figure \ref{fig:LDE} shows the linear dispersion error for $N_x = 2$ for an increasing $kh$.
Unsurprisingly, it is clear that the waves are well-resolved with only a few spatial modes in the horizontal when using a trigonometric basis in the horizontal direction, as all wave numbers exceeding the Nyquist frequency for the shortest wavelengths that can be resolved with as little as two grid points ($N_x =2$) in such a Fourier basis. The model is also capable of simulating cases of waves corresponding to a high dispersion $kh$, i.e. short wavelength compared to depth, with relatively few nodes. If we compare the results to spectral accurate potential flow models such as Engsig-Karup et al. \cite{Allan16}, we see that the model performs excellently in comparison, noting that our model has a higher range of $kh$ with low dispersion error for the same vertical resolution in comparison to the full high order model in Engsig-Karup et al., and avoids the dispersion errors at low $kh$ as seen in the mixed-order model. 

\subsection{Oscillatory wave boundary layer}
\begin{figure}[H]
\vspace*{0cm}
\hspace*{0cm}
\centering
\includegraphics[width=110mm]{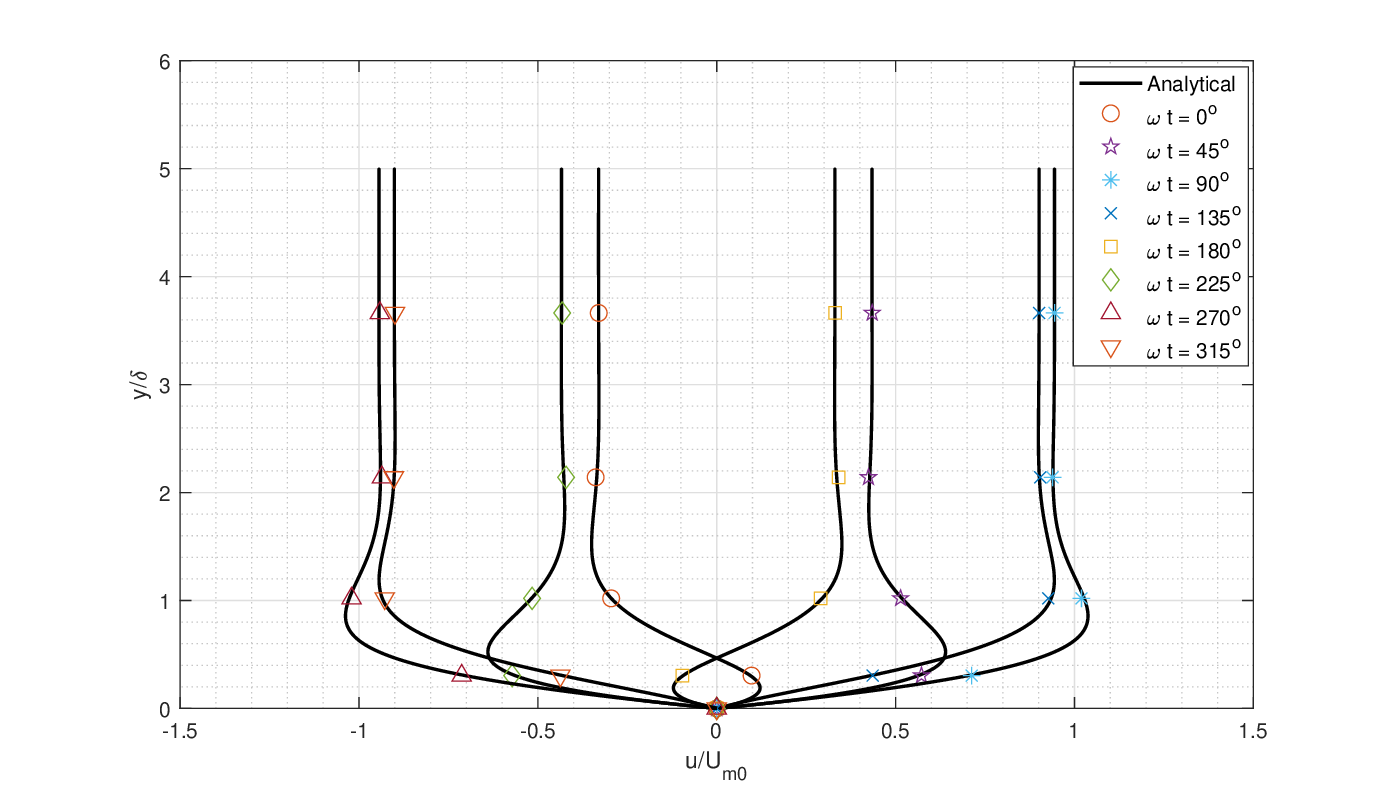}
\caption{Wave boundary velocity profile}
\label{fig:LWB}
\end{figure}

A main characteristic of the NSE solver that differs from shallow water and potential flow solvers is the ability to describe the boundary layer taking into account viscous effects, which is indispensable for the case of sediment transport and drag and friction analysis of offshore structures. As we will show below, the proposed model is capable of resolving the boundary layer. A no-slip boundary condition is imposed at the sea bed 
\begin{align}
    \textbf{u} = 0\quad \text{on} \quad \Gamma_b,
\end{align}
with $\Gamma_b$ being the solid bottom boundary.
The analytical solution of the oscillatory wave boundary layer can be derived from the small-amplitude wave theory \cite{Sumer20}:
\begin{subequations}
\begin{align}
    u(z,t)=U_{0m}\mbox{cos}(\omega t-k x)-U_{0m} \exp\left(-\frac{z}{\delta_1}\right)\mbox{cos}\left(\omega t-k x-\frac{z}{\delta_1}\right),
\end{align}
where 
\begin{align}
    U_0&=U_{0m}\mbox{cos}(\omega t-k x),\\
    U_{0m}&=\frac{\pi H}{T}\frac{\mbox{cosh}(z+h)}{\mbox{sinh}(k h)},\\
    \delta_1&=\left(\frac{2\nu}{\omega}\right)^{0.5}.
\end{align}    
\end{subequations}
Here $U_0$ is the horizontal part of the orbital velocity just outside the boundary layer, $U_{0m}$ is the amplitude of the orbital velocity and $\delta_1$ is the Stokes length. 
The numerical results are computed with the linearised NSE model as defined in \eqref{eq:NSElin} and \eqref{eq:linFS} with $N_x = 20$ and $N_\sigma=50$. The wave is initialized from linear wave theory with $H=0.02$ m and $kh = 0.6725$ on a periodic grid. 
The predicted velocity profile by the model can be seen in Figure \ref{fig:LWB}, where it is compared against the analytical solution. Note that $\delta = \frac{3\pi}{4}\delta_1$ defines the boundary layer thickness. As can be seen, the model is able to resolve the boundary layer despite using very few grid points within the boundary layer itself. 

\subsection{Harmonic generation over a submerged bar}
\begin{figure}[H]
    \centering
    \includegraphics[width=120mm]{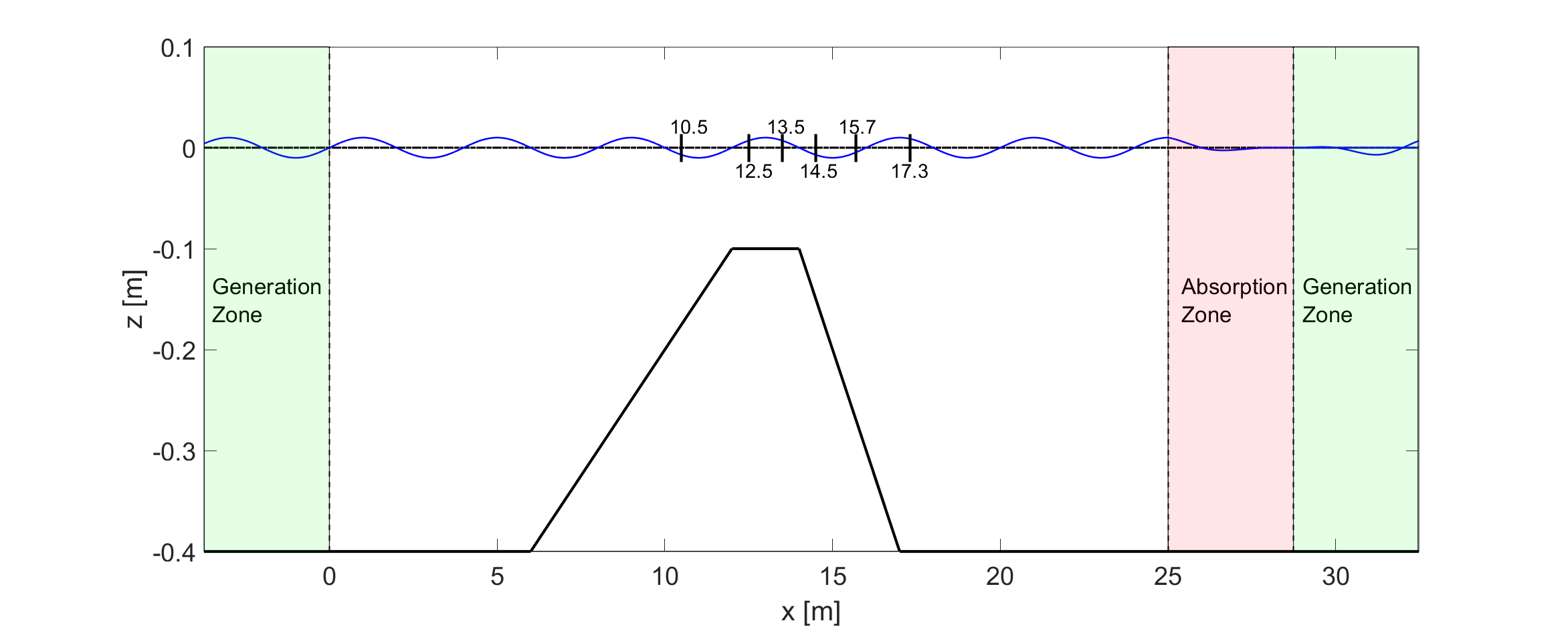}
    \caption{Submerged bar wavetank setup. The numbered bars denote the measurement location.}
    \label{fig:barsub}
\end{figure}

\begin{figure}[H]
\vspace*{0cm}
\hspace*{0cm}
\centering
\subfloat[]{
  \includegraphics[width=50mm]{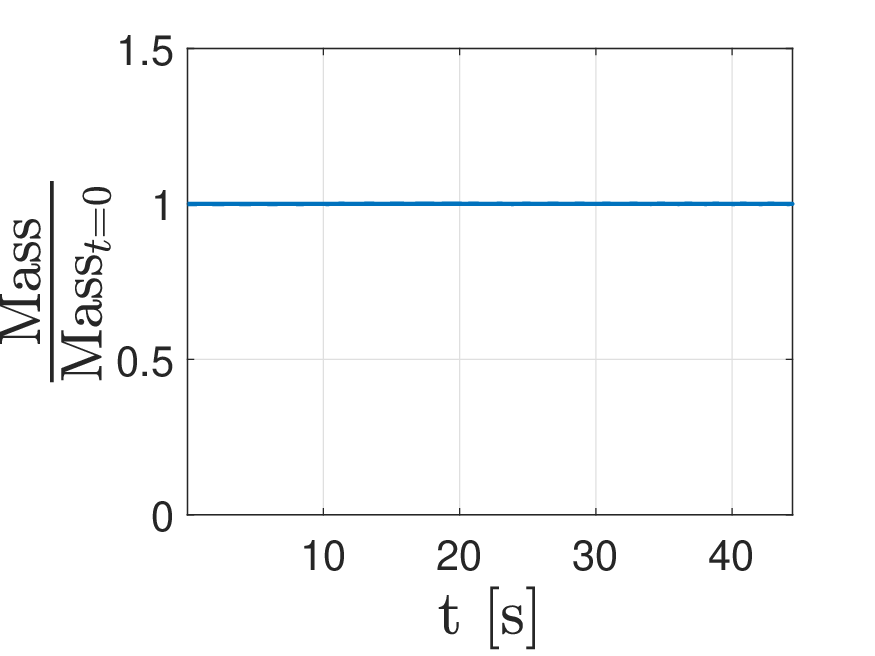}
}
\subfloat[]{
  \includegraphics[width=50mm]{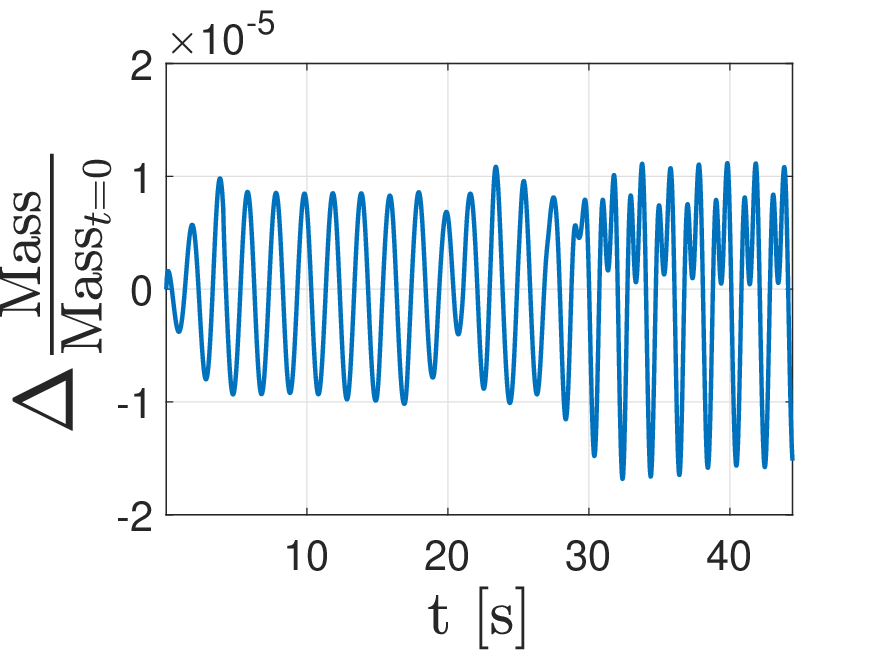}
}
\caption{Time series of the overall mass of the system (a) and rate of change in mass (b).}
\label{fig:barmass}
\end{figure}

Harmonic generation over a submerged bar is a well-known test problem, first posed by Beji and Battjes \cite{Beji93,Beji94}, who performed both experimental and numerical experiments. The setup is a wave-tank with a non-flat bottom, which can be seen in Figure \ref{fig:barsub}. 
The waves are generated through a relaxation method due to \cite{Larsen83} and is used for setting up a generation and absorption zones in line with the work due to Engsig-Karup \cite{Allanphd}. Note that the generation zone has been extended to the right of the domain to ensure periodicity required due to the use of FFT. The zones each have the same length as the initial waves. The initially generated waves are based on the stream function, and have a dimensionless depth of $kh = 0.67$, with an initial wave height of $H = 0.02$ m, and a wave period of $T = 2.02$ s, leading to an initial wavelength of $L = 3.737$ m. The waves incur a shoaling effect during the propagation over the bar, which causes the waves to steepen. After the top of the bar, the waves then decompose into shorter free waves with rapid changes to the wave structure. 
\begin{figure}[H]
\vspace*{0cm}
\hspace*{0cm}
\centering
\subfloat[]{
  \includegraphics[width=50mm]{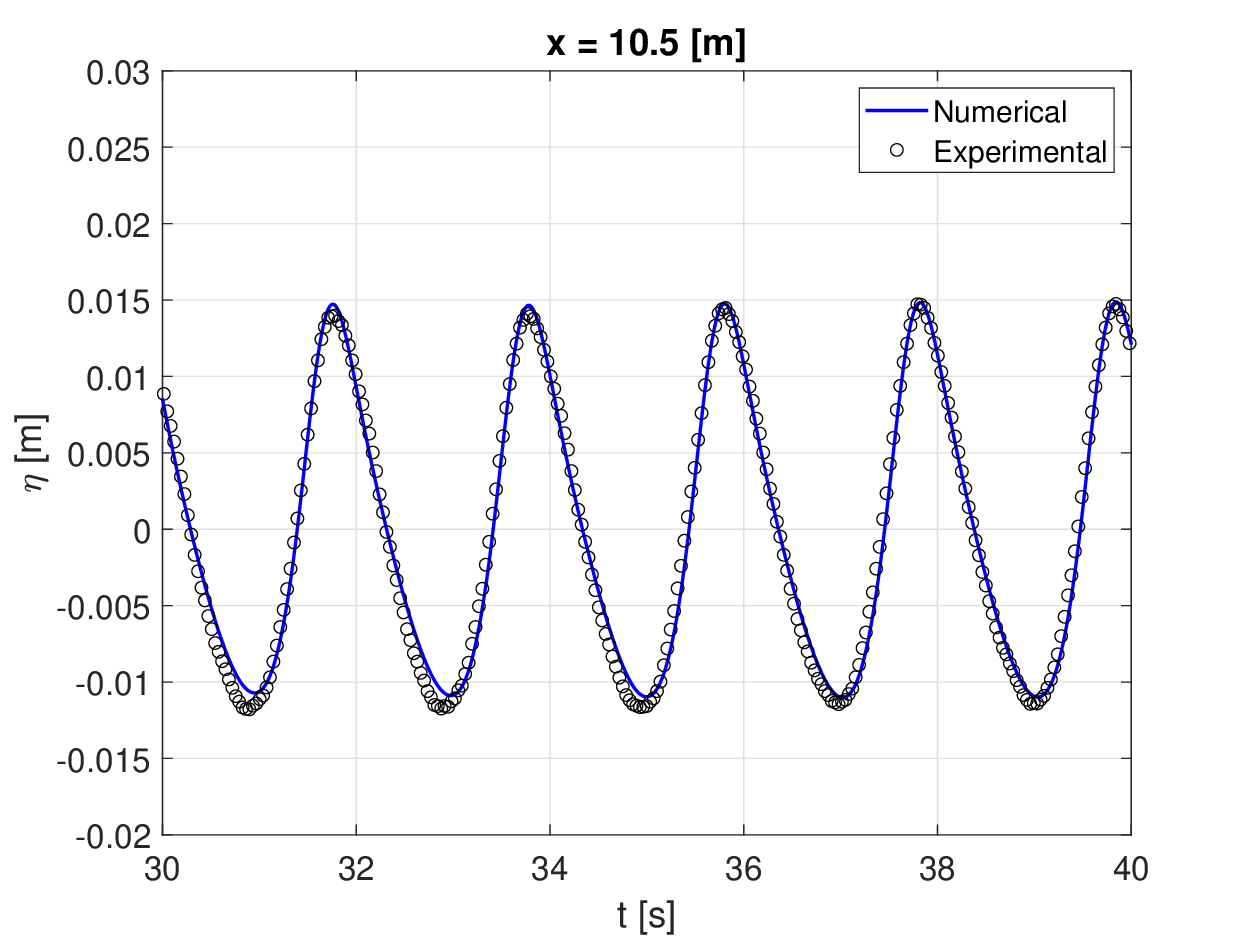}
}
\subfloat[]{
  \includegraphics[width=50mm]{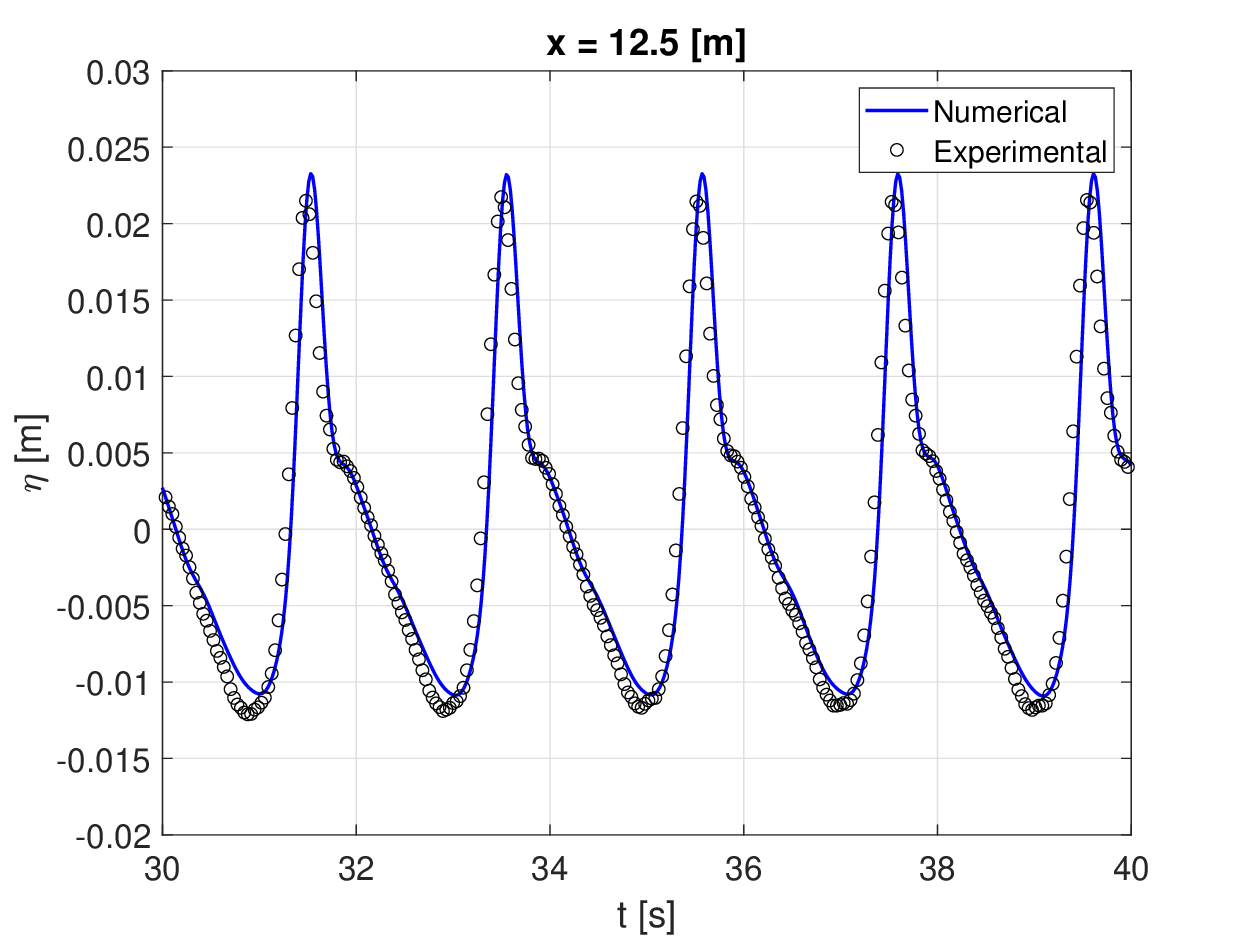}
}
\hspace{0mm}
\hspace*{0cm}
\subfloat[]{
  \includegraphics[width=50mm]{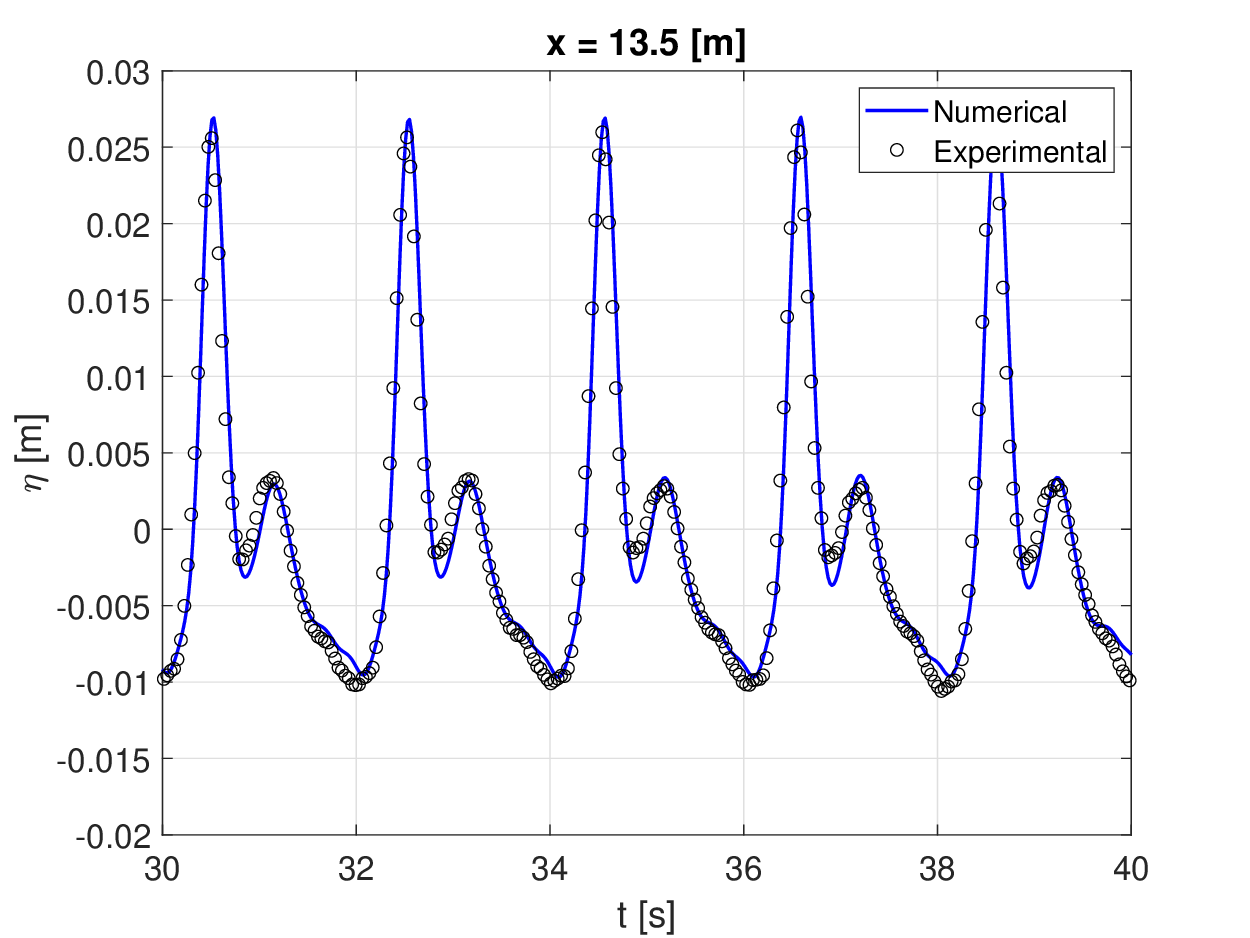}
}
\subfloat[]{
  \includegraphics[width=50mm]{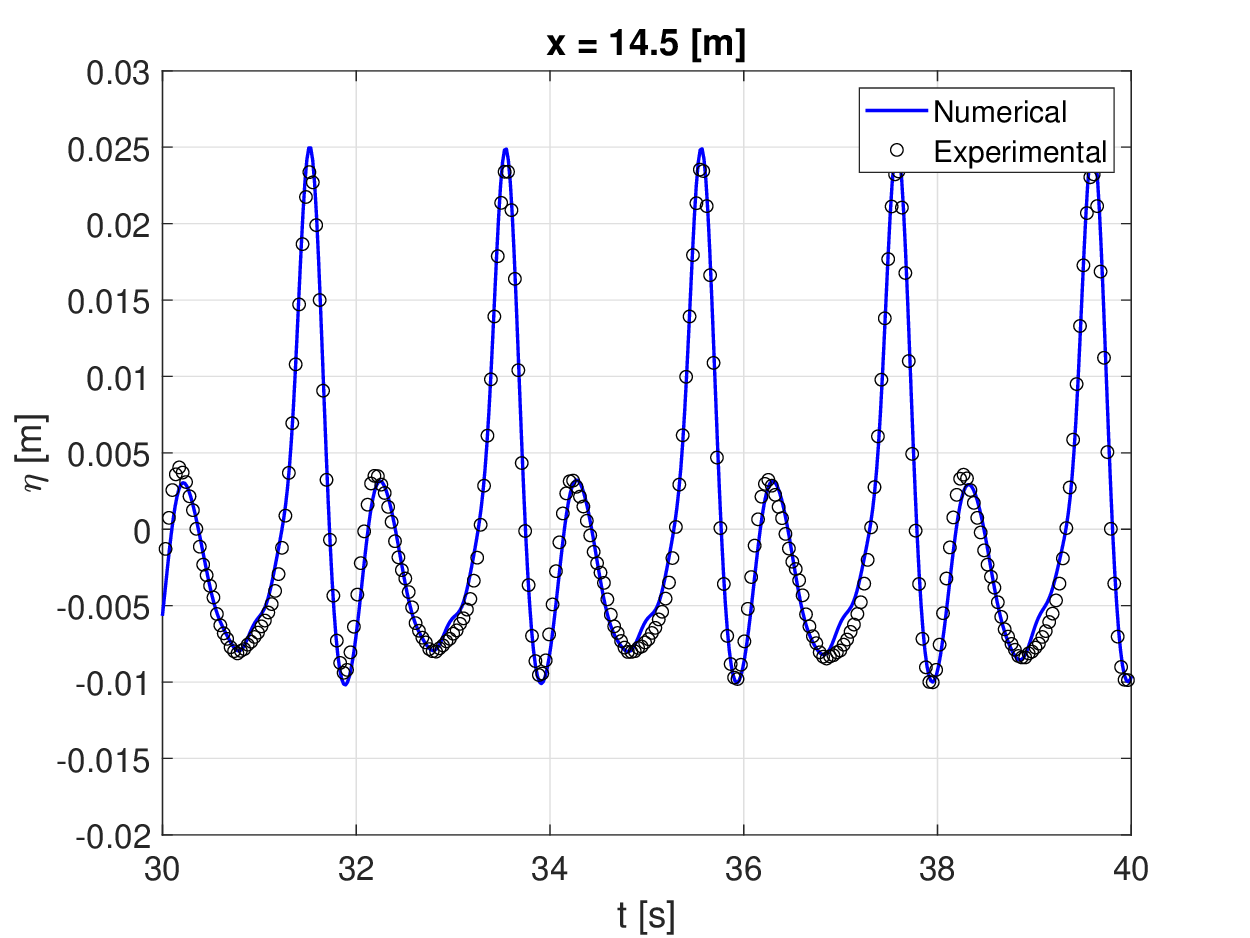}
}
\hspace{0mm}
\hspace*{0cm}
\subfloat[]{
  \includegraphics[width=50mm]{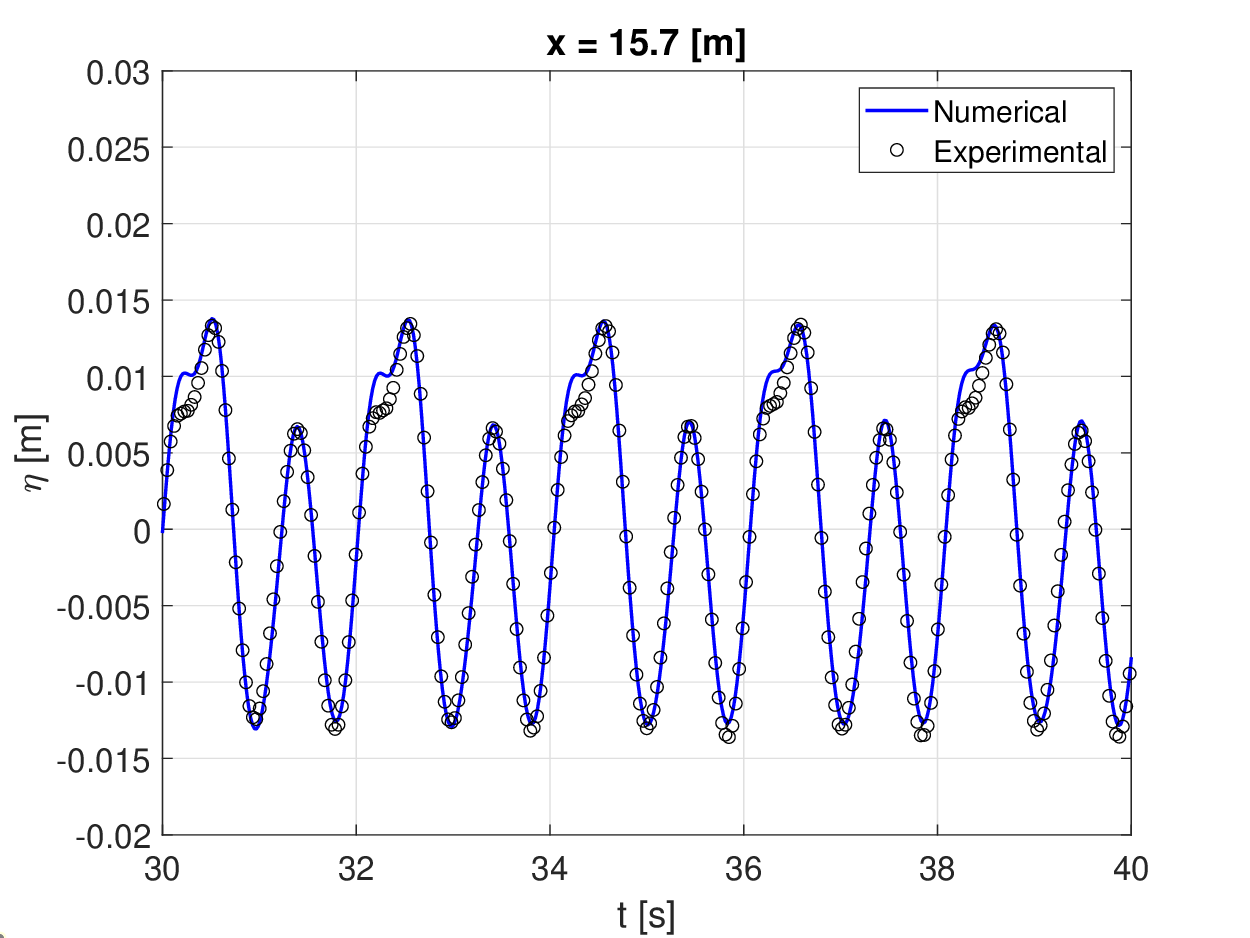}
}
\subfloat[]{
  \includegraphics[width=50mm]{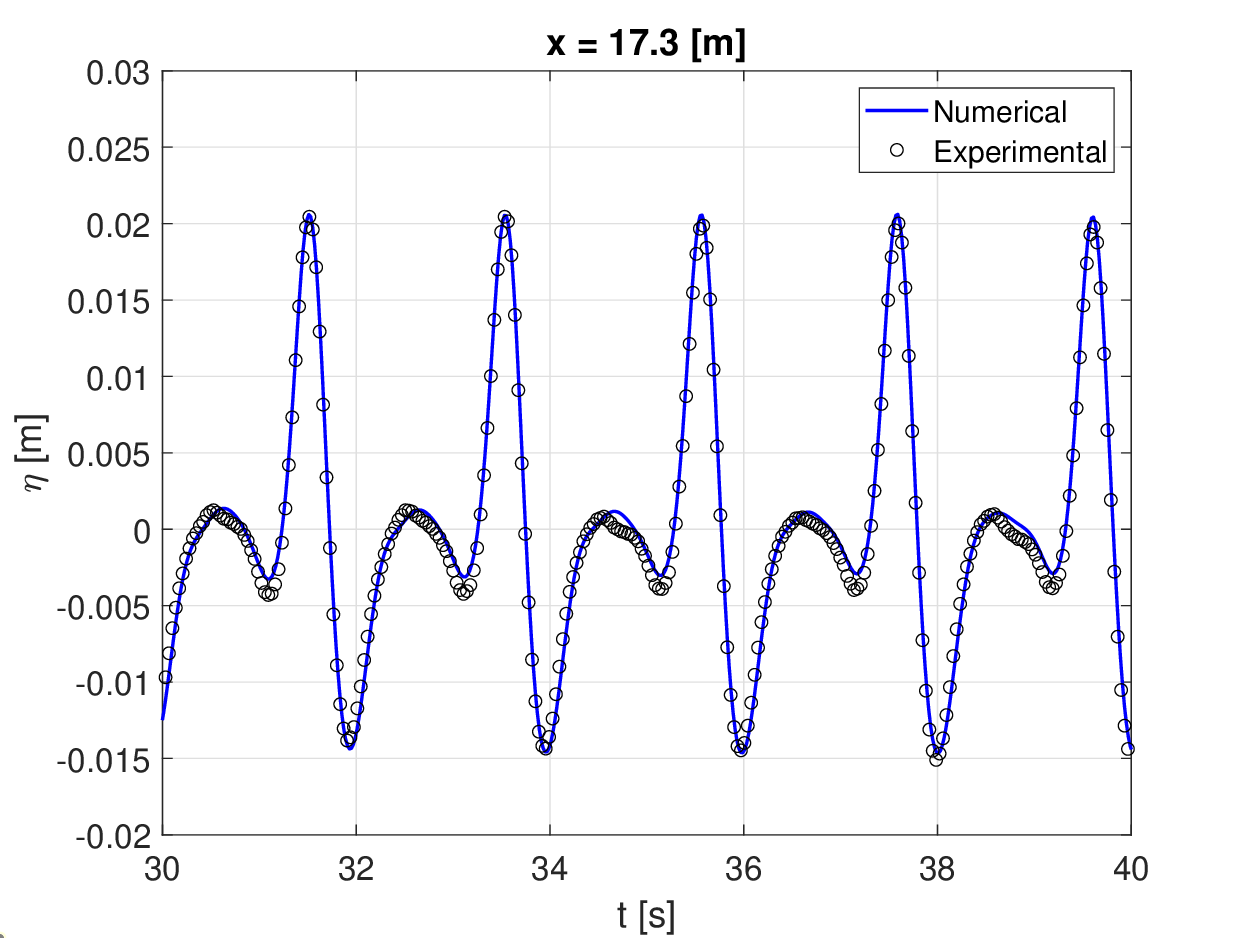}
}

\caption{Time series of the free surface elevation at measuring points for both numerical and experimental results.}
\label{fig:barres}
\end{figure}
As shown in Figure \ref{fig:barres}, the comparison to experimental results shows that the solver is capable of capturing the evolution of the nonlinear and dispersive waves accurately. Also, the results are found to compare favorably to Boussinesq models due to Gobbi \& Kirby \cite{Gobbi98}, and to those of nonlinear potential flow models such as Engsig-Karup et al. \cite{Allan16}. It is noted that the viscous and rotational effects accounted for in the incompressible Navier-Stokes model used in this study does not significantly differ from results compared to these other works. Moreover, as shown in Figure \ref{fig:barmass}, the solver maintains mass conservation throughout the simulation. We note that the rate of change shows a wave pattern, which is due to the generation and absorption zones being slightly offset in terms of adding and removing mass.

\section{Conclusion}
We present a new nodal hybrid-spectral method that solves the incompressible Navier-Stokes problem with a free surface. The method achieves the desired spectral convergence for highly-nonlinear waves at both shallow, deep, and intermediate water and achieves very low numerical dispersion error. The method is also capable of propagating waves over uneven bottom, with results matching those of experiments, and can resolve boundary layers with very few grid points. Moreover, it is shown that utilizing a geometric $p$-multigrid method as preconditioning for the GMRES method results in significant reductions in both iteration count for the iterative solver and computation time. In ongoing work, the new scheme will investigate turbulent flows and also be prepared for large-scale sea-state estimation in regional areas. A limitation of the hybrid-spectral scheme proposed in this work is that it is not straightforward to handle floating offshore structures and is is therefore a topic of ongoing work to develop a new high-order spectral element scheme that provides a basis for such advanced modelling of wave-structure interaction problems.

\section{Acknowledgments}
The research was funded and hosted by the Department of Applied Mathematics and Computer Science (DTU Compute). We thank the DTU Computing Center for access to HPC resources.





\bibliography{sn-bibliography}


\end{document}